\documentclass[11pt]{article}

\usepackage{amsmath, amssymb, amscd, amsthm, mathrsfs, url, textcomp, bbm}
 \usepackage{amssymb}
\usepackage{graphicx}
\usepackage{enumitem}
\usepackage{cleveref}
\usepackage{leftidx}
\usepackage{color}

\allowdisplaybreaks

\usepackage{ifxetex}
\ifxetex
  \usepackage{fontspec}
\else
  \usepackage[T1]{fontenc}
  \usepackage[utf8]{inputenc}
  \usepackage{lmodern}
\fi

\DeclareMathAlphabet{\mathpzc}{OT1}{pzc}{m}{it}

\usepackage{amsfonts}
\usepackage{geometry}
\usepackage{dsfont}

\numberwithin{equation}{section}

\newcommand{\be}{\begin{eqnarray}}
\newcommand{\ee}{\end{eqnarray}}
\newcommand{\ce}{\begin{eqnarray*}}
\newcommand{\de}{\end{eqnarray*}}
\newtheorem{theorem}{Theorem}[section]
\newtheorem{lemma}[theorem]{Lemma}
\newtheorem{remark}[theorem]{Remark}
\newtheorem{definition}[theorem]{Definition}
\newtheorem{proposition}[theorem]{Proposition}

\newtheorem{corollary}[theorem]{Corollary}

\def\1{{\mathbbm 1}}

\newcommand{\bP}{\mathbb{P}}
\newcommand{\fP}{\mathbf{P}}

\newcommand{\bE}{\mathbb{E}}
\newcommand{\fE}{\mathbf{E}}

\def\wt{\widetilde} 

\newcommand{\R}{\mathbb{R}}

\newcommand{\N}{\mathbb{N}}

\newcommand{\sL}{\mathcal{L}}
\newcommand{\sM}{\mathcal{M}}
\newcommand{\sN}{\mathcal{N}}

\def\eps{\varepsilon}

\title{Harnack inequality for weakly coupled nonlocal systems}
\author {Zhen-Qing Chen \quad  Xiangqian Meng}
\date{\today} 

\providecommand{\keywords}[1]
{
  \textbf{Keywords: } #1
}

\begin{document}
\maketitle
\begin{abstract}
In this paper, we consider a weakly coupled system of nonlocal operators which contain both
diffusion part with uniformly elliptic diffusion matrices and bounded drift vectors and the jump part
with relatively general jump kernels. We use the two-sided scale-invariant Green function estimation to
prove the scale-invariant Harnack inequality for the weakly coupled nonlocal systems. In the case where the switching rate matrix is
strictly irreducible, the scale-invariant full rank Harnack inequality is proved.Our approach is mainly probabilistic.
\end{abstract}

\keywords{Weakly coupled nonlocal system; Harnack inequality; H\"older regularity; Nonlocal operators; Green function estimates; conditional processes}

{\small
 \begin{tableofcontents}
 \end{tableofcontents} }
\section {Introduction}
 Let $\sM=\{1,...,m\},\hbox{ for } m\geq 2$. Suppose that $Q(x)=(q_{ij})_{m\times m}(x)$ is an $m\times m$ matrix valued function on $\R^d$ such that
for a.e. $x\in \R^d$, 

$$q_{ij}(x)\geq  0 \hbox{ for } i\neq j, \quad\hbox{ and } \sum_{j=1}^m q_{ij}(x)\leq 0 \hbox{ for each }i\in \sM. $$

For any function 
\begin{equation*}
    u = \left(\begin{matrix}
     u_1\\
    \vdots\\
     u_m
     \end{matrix}
     \right):\R^d\mapsto \R^m, 
 \end{equation*}
 define 
\begin{equation}\label{e:1.1}
Q(x)u(x,\cdot)(i):=\sum_{j\in \sM}q_{ij}(x)u(x,j),\quad (x,i)\in \R^d\times \sM.
\end{equation} 

  Consider the following weakly coupled operator $\mathcal{G}$ such that for each $i\in \sM, u(\cdot, i)\in C_b^2(\R^d)$, 
\begin{equation}\label{e:1.2}
\mathcal{G}u(x,i):=\sL_iu(x, i)+Q(x)u(x,\cdot)(i),\quad (x,i)\in \R^d\times \sM
\end{equation}
with
\begin{eqnarray}\label{e:1.3}
\sL_iu(x,i)&=&\sL_i^{b_1}u(x,i)+\mathcal{S}^{b_2}u(x,i),
\end{eqnarray}
where
\begin{eqnarray}\label{e:1.4}
\sL_i^{b_1}u(x,i):&=&\frac{1}{2}\sum_{k,l=1}^da_{kl}(x, i)\frac{\partial^2}{\partial x_k\partial x_l}u(x,i)+\sum_{k=1}^db_1^k(x,i)\frac{\partial}{\partial x_k}u(x,i);\\
\mathcal{S}^{b_2}u(x,i):&=&\int_{\R^d\setminus \{0\}}(u(x+z,i)-u(x,i)-\nabla u(x,i)\cdot\mathds{1}_{|z|\leq 1}z)b_{2}(x,z,i)j_i(z)dz.
\end{eqnarray}
Here $a_{kl}(x,i)=a_{lk}(x,i)$ and there are positive constants $\Theta_1,\Theta_2, \Theta_3,\Theta_4,c,c_1$ and $\gamma\in (0,1)$  such that  $(a_{kl}(x,i))_{1\leq k,l\leq m},(b^k_1(x,i))_{1\leq k\leq m}, b_2(x,z,i),j_i(z)$ satisfy the following conditions:

\begin{description}
\item{(a)} (Uniform ellipticity and H\"older continuity on diffusion matrix) 
\begin{eqnarray}\label{e:2.4}
\Theta_1 |\zeta|^2\leq \sum_{k,l=1}^{d} a_{kl}(x,i) \zeta_k \zeta_l&\leq&\Theta_1 ^{-1} |\zeta |^2, \quad \hbox{ for } \zeta  \in \R^d;\\
|a_{kl}(x_1,i)-a_{kl}(x_2,i)|&\leq& c|x_1-x_2|^{\gamma}, \hbox{ for }(x_1,i),(x_2,i) \in \R^d \times \sM,\hbox{ some }\gamma\in (0,1).\nonumber\\
\end{eqnarray}
\item{(b)}  (Bounded drift)
\begin{equation}\label{e:2.5}
\|b_1 \|_{\infty}:=\sup_{(x,i)\in \R^d\times \sM} |b_1(x,i)| \leq \Theta_2.
\end{equation}
\item{(c)} (L\'evy jumping kernel condition and nonnegative local boundedness) 
\begin{eqnarray}\label{e:2.5a}
\|b_2(x, z,i)\|_{L^{\infty}(\R^d\times \R^d\setminus\{0\}\times \sM)}\leq \Theta_3.
\end{eqnarray}
\begin{eqnarray}\label{e:2.5b}
\int_{\R^d\setminus\{0\}}(1\wedge |z|^2)j_i(z)dz<\Theta_4, 
\end{eqnarray}
where each $j_i(z)$ is a nonnegative locally bounded function on $\R^d\setminus \{0\}$ such that there exists a constant $\beta\in (1,2)$ and $c_1>0$ satisfying
\begin{eqnarray}\label{e:2.6}
j_i(z)&\leq& \frac{c_1}{|z|^{d+\beta}}  \hbox{ for }z\in B(0,1)\setminus \{0\} , i\in \sM.
 \end{eqnarray}
\end{description}
In matrix form, 
\begin{equation}\label{e:1.5}
 \mathcal{G}u=
 \left(\begin{matrix}   \sL_{1} & 0 & \dots & 0 \\0 & \sL_2 & \dots & 
0\\ \vdots & \vdots & \ddots & \vdots\\0 & 0 & \cdots & \sL_m \end{matrix}\right)u+Qu.
\end{equation}
Such operators with no derivatives appearing in the coupling terms are called  {\it weakly coupled } operators. There are studies on potential theory for weakly coupled elliptic systems where $\sL_i$ in \eqref{e:1.3} are all  differential operators. In 1989, Skorohord in \cite{skorohod} introduced a Markov process $\{(X_t, \Lambda_t),t\geq 0; \bP^{(x, i)}, 
(x,i) \in (\R^d\times \sM)\}$in a descriptive manner under the condition that $\sum_{j\in \sM}q_{ij}(x)=0 $ a.e. in  $\R^d$ having $\mathcal{G}$ as its infinitesimal generator \eqref{e:1.2}. In \cite{switch}, Chen and Zhao established the existence of switched diffusion processes corresponding to the weakly coupled elliptic operators in divergence form with measurable coefficients such that the diffusion matrix is uniformly elliptic, and the coefficients of the first order term and each term in the switching rate matrix $Q$ belongs to a Kato class. They also showed the solvability of the Cauchy problems for weakly coupled elliptic systems via a Dirichlet space approach, i.e. $\partial u/\partial t=\mathcal{G}u$ and a strong positivity result for the solutions. In \cite{switch}, the regularity condition on coefficient $a^k$, $b^k$, and $q_{kl}$ are the following :
\begin{eqnarray*}
&&a_{kl}(\cdot, i)\in W^{1,2}(D)\hbox{ for each i}, \sum_{k=1}^d(\partial b_k(\cdot, i))/\partial x_k)\hbox{ is bounded from below for each }i, \\
&&\hskip 0.5truein  \hbox{ and }\sum_{i=1}^mq_{ij} \hbox{ is bounded from above for each }j,
\end{eqnarray*}
 where $D$ is a bounded domain. In \cite{potentialsystem}, Chen and Zhao generalized the previous result to weakly coupled elliptic system (i.e. $\mathcal{G}u=0$) whose coefficients are only measurable such that the diffusion matrix is uniformly elliptic with  lower-order terms in a Kato class on bounded domains. They also gave a probabilistic representation of the solution and proved a strong positivity result for that problem via a probabilistic approach. Later, the corresponding Harnack inequality for weakly coupled elliptic systems and the full Harnack inequality under the irreducibility assumption of $Q$ and  H\"older continuous coefficients are established by an analytic method in \cite{weaklycoupled}. The Harnack inequality improves the result of F. Mandras \cite{mandras}, in which the Harnack inequality takes the form :
$$\sum_{i=1}^m \sup_{x\in D} u_i(x)\leq C\sum_{i=1}^m \inf_{x\in D} u_i(x),$$ where $D$ is bounded. The approach in \cite{weaklycoupled} is based on the representation theorem obtained in \cite{potentialsystem} and the estimates of Green functions and harmonic measures of the operators $\sL_k+q_{kk}$ in small balls. Stimulated by the ongoing research on nonlocal operators which correspond to various jump processes, see \cite{Bass,basslevin,caffarelli,CK03,renminsong}, as well as the surgent needs of applications of switched Markov processes $(X, \Lambda)$ to the fields of engeering and finance \cite{stock, meshfree, mari, algorithm}, the scale-invariant Harnack inequality for $\sL_i$ of the form \eqref{e:1.3}  with uniformly elliptic diffusion matrix, bounded drift and more general jump kernel assumption has been studied in \cite{hinonlocal}, and the weakly coupled nonlocal system of \eqref{e:1.2} has also been studied, for example, in \cite{jumpswitching}. In \cite{jumpswitching}, the non scale-invariant Harnack inequality for each level and the full rank Harnack inequality for the operator  \eqref{e:1.2}  satisfying \eqref{e:2.4},\eqref{e:2.5},\eqref{e:2.5b}, the additional continuity assumption on $(a_{kl}(x))$ and $(b_k(x))$, and the relative comparison of the jump intensity have been established. The purpose of this paper is to establish the scale-invariant Harnack inequality for the weakly coupled system.  In this paper we assume that the dimension $d\geq 2$ and our goal is to establish the scale-invariant Harnack inequality for harmonic functions with respect to the operator $\mathcal{G}$ in \eqref{e:1.2}.
  We say that the non-local operator $\sL$ of the form \eqref{e:1.3} belongs to {\bf class $\sN(\Theta_1, \Theta_2, \Theta_3,\Theta_4,\gamma, \beta, c_1)$}
if the following \eqref{e:2.4}-\eqref{e:2.6} hold and an operator $ \mathcal{G}$ belongs to {\bf class $\sN(c_0, m,\gamma, \beta, \Theta_1, \Theta_2, \Theta_3, \Theta_4,\Theta_5, c_1, \vartheta)$ }if \eqref{e:2.4}-\eqref{e:2.6} hold with the additional conditions \eqref{condQ1} and \eqref{UJScondition} for  the switching rate matrix $Q$:
\begin{eqnarray}\label{conditiononQ}
&&\hbox{There exists a constant matrix }Q^0=(q_{ij}^0) \hbox{ and a positive constant }0<c_0<1\nonumber\\
&& \hskip 0.5  truein  \hbox { such that } c_0q^0_{ij}\leq q_{ij}(x)\leq q^0_{ij}\leq \Theta_5 ,i\neq j   \, \, a.e. \, \hbox{ on } \R^d \label{condQ}\\
&&  \hskip 0.5  truein \hbox { and } -q_{ii}(x)\leq -q_{ii}^0 \leq\Theta_5 \, \, a.e. \, \hbox{ on } \R^d \label{condQ1}.
\end{eqnarray}

Define $J^{b_2}(x, y,i)=b_{2}(x,y-x,i)j_i(y-x),$ for each $x,y\in \R^d, i\in \sM$. 
We say that the jumping kernel $J^{b_2}(x, x+z,i),i\in \sM$ satisfies the {\bf UJS condition} if there exists a constant $\vartheta=\vartheta(d, \beta, b_2)>0$ so that for a.e. $x\in \R^d,z \neq 0$, every $i\in \sM$, 
\begin{equation}\label{UJScondition}
J^{b_2}(x,x+z,i)\leq \frac{\vartheta}{|B(x,r)|}\int_{B(x,r)}J^{b_2}(u,{x+z},i)du \hbox{ for every } 0<r<|z|/2.
\end{equation}
Here $|B(x_0, r)|$ denotes the Lebesgue measure of the ball $B(x_0,r)$, which is $\omega_dr^d$ with $\omega_d$ being the volume of the unit ball in $\R^d$. For the typical examples of the nonlocal operators whose jump kernel $J^{b_2}$ satisfy the UJS condition, we refer the readers to Remark 1.6 in \cite{BHP}.

Given a switching rate matrix $Q$, we denote the collection of $n$-step path from state $i$ to $j$ by
 \begin{eqnarray}\label{path1}
 \Psi(n;i,j):&=&\{(l_0,l_1,...,l_n):l_0=i,l_n=j,l_i\in \sM, l_i\neq l_{i+1},\{q_{l_il_{i+1}}(x)>0\} \nonumber\\
 &&\hbox{ has a positive Lebsgue measure on }\R^d, \hbox{ for }i=0,1,...,n-1\}.
\end{eqnarray}
By \eqref{conditiononQ},$\Psi^0(n;i,j)$, which is defined correspondingly  for $Q^0$, satisfies that $\Psi^0(n;i,j)=\Psi(n;i,j)$ for any $i,j\in \sM, n\in \N$.
\begin{definition}
We say that the operator $\mathcal{G}$ or its associated matrix-valued function $Q$ is {\bf{irreducible}}  in $\R^d$ if for any distinguished $k,l\in\sM$, there exist  $l_0,...,l_n$ in $\mathcal{M}$ with $l_{i-1}\neq l_i$ for $1\leq i\leq n$, $l_0=k$ and $l_n=l$ such that $q_{l_{i-1}l_i}^0\neq 0.$
\end{definition}

\begin{definition}
We say that the operator $\mathcal{G}$ or its associated matrix-valued function $Q$ is {\bf{strictly irreducible}} if $q^0_{kl}>0$ for any $k,l\in \sM$. 
\end{definition}
And we denote
\begin{equation}\label{lowerbdofQ0}
q_0=\min\{q^0_{ij}: i,j\in \sM,i\neq j \}.
\end{equation}

For any $\mathcal{G} \in \sN(c_0, m,\gamma, \beta, \Theta_1, \Theta_2, \Theta_3, \Theta_4,\Theta_5, c_1, \vartheta)$, we can show that there is a unique Hunt process $Y=\{(X_t,\Lambda_t), t\geq 0;\bP^{(x,i)},(x,i)\in (\R^d\times \sM)\}$ solving the martingale problem for $(\mathcal{G},C^2_b(\R^d\times \mathcal{M}))$ starting at $(x,i)$ :  that is, for every function $u(x, i):\R^d\times \mathcal{M}\mapsto \R$, such that $u(\cdot, i)\in C_b^2(\R^d)$ for every $i\in \sM$, 
$$M_t^u:=u(X_t,\Lambda_t)-u(X_0,\Lambda_0)-\int_0^t\mathcal{G}u(X_s,\Lambda_s)ds$$
is a $\bP^{(x,i)}$-local maringale for each $(x,i)\in \R^d\times \sM$ and $\bP^{(x,i)}(X_0=x, \Lambda_0=i)=1$, where $X$ is an $\R^d-$valued process, and $\Lambda$ is a switching process taking values  in $\mathcal{M}$. In this case, we say that the Hunt process $Y=\{(X_t,\Lambda_t), t\geq 0;\bP^{(x,i)},(x,i)\in (\R^d\times \sM)\}$ {\bf has the generator $\mathcal{G}$ or corresponds to the operator} $\mathcal{G}$.

 The process $Y$ can be constructed in the way described in \cite[Remark 2.2]{weaklycoupled}: Under the assumption \eqref{e:2.4}-\eqref{e:2.6}, by Theorem 5.2  in \cite{komatsu}, there is a unique conservative strong Markov process $\bar{X}=\{\bar{X}_t, t\geq 0;\bP^x, x\in \R^d\}$ whose distribution with $\bP^x(X_0=0)=1$ solves the martingale problem for $(\sL, C^2_b(\R^d))$, where $\sL\in \sN(\Theta_1, \Theta_2, \Theta_3,\Theta_4, \gamma, \beta, c_1)$. By the assumption \eqref{e:2.5a}, we know that $b(x,z)j(z)dz\leq \Theta_3j(z)dz$ for every $x\in \R^d$, then by \cite[Theorem 2.1, Theorem 2.2]{komatsu}, the process $\bar{X}$ is quasi-left continuous under $\bP^x$. Fix $i_0\in \sM$ and let $(X^{i_0}, \bP^{(x_0,i_0)})$ be the Hunt process corresponding  to  the infinitesimal generator $\sL_{i_0}$ starting from $(x_0, i_0)$.We run the subprocess $\wt{X}^{i_0}$ of $X^{i_0}$ that is killed with the rate $-q_{i_0i_0}(x)$. Note that the subprocess $\wt{X}^{i_0}$ has the infinitesimal generator $\sL_{i_0}+q_{i_0i_0}$. At the lifetime $\tau_1$, the killed subprocess $\wt{X}^{i_0}$  is killed with probability $1+\sum_{j\in \sM}q_{i_0j}(X^{i_0}_{\tau_1-})/q_{i_0i_0}(X^{i_0}_{\tau_1-})$ and jumps to plane $j\neq i$ with probability $-q_{i_0j}(X^{i_0}_{\tau_1-})/q
_{i_0i_0}(X^{i_0}_{\tau_1-})$ and then starting from $X^{i_0}(\tau_1-)$, we run an independent copy of a subprocess $\wt{X}^{j}$ of $X^j$ with the killing rate $-q_{jj}(x)$. Iterating this procedure, the resulting process $Y=((X_t, \Lambda_t),t\geq 0;\bP^{(x,i_0)})$ is a Hunt proces with lifetime $\zeta$ by \cite{ikeda, merge}. For each $x\in \R^d$, we say that $Q(x)$ is {\it Markovian} if $\sum_{j\in \sM} q_{ij}(x)=0$ a.e. on $\R^d$ for every $i\in \sM$, and {\it sub-Markovian}  if $\sum_{j\in \sM }q_{ij}(x)\leq 0$ a.e. on $\R^d$ for every $i\in \sM$. When $Q(x)$ is {\it Markovian}, the lifetime $\zeta=\infty$, and when $Q(x)$ is {\it sub-Markovian}, $\zeta<\infty$ with positive probability. We use the convention $(X_t, \Lambda_t)=\partial$ for $t\geq \zeta$, where $\partial$ is a cemetery point and any function is extended to $\partial$ by taking value zero. It is easy to check that $((X, \Lambda), \bP^{(x_0, i_0)})$ solves the maringale problem for $(\mathcal{G}, C_b^2(\R^d\times \sM))$, where we denote by $C_b^2(\R^d\times \sM)$  the class of functions $u$ defined on $\R^d\times \sM$ such that $u(\cdot, i)$ is in $C_b^{2}(\R^d)$ for each $i\in \sM$. This way of constructing a switched diffusion process by patching together the pre-switching process $\wt{X}^{i}$ for the operator $\sL_i+q_{ii}$ with its switching distribution has also been utilized in \cite[Page 296]{potentialsystem}.

In this paper, we always assume that $D$ is a bounded and connected open set in $\R^d$.  We say that  a Borel measurable function $u:D\times \sM\mapsto \R^d$ is {\textbf{$\mathcal{G}$-harmonic}} in $D\times \sM$ if for any relatively compact open subset $V$ of $D\times \sM$, $$\bE^{(x,i)}[u(X_{\tau_V}, \Lambda_{\tau_V})]<\infty, \hbox{ and } u(x,i)=\bE^{(x,i)}[u(X_{\tau_V}, \Lambda_{\tau_V})]$$ for every $(x,i)\in V$,
where $\tau_V:=\inf\{t\geq 0:(X_t, \Lambda_t)\notin V\}$ is the first exit time from the set $V$. Heuristically, $u$ is $\mathcal{G}$-harmonic in $V$ if $\mathcal{G} u=0$ in $V$. But we are not going to establish this analytic characterization.  See \cite{harmonic, mazhu} for the equivalent characterizations between probabilistic and analytic notions of harmonicity under some suitable conditions.

Our first main result is on the local H\"older regularity for bounded $\mathcal{G}$-harmonic functions.
 \begin{theorem}\label{T:Holder} Let $\mathcal{G}\in \sN(c_0, m,\gamma, \beta, \Theta_1, \Theta_2, \Theta_3, \Theta_4,\Theta_5, c_1, \vartheta)$.
 There exist constants  \hfill \break
  $\wt{r}_0=\wt{r}_0(d, \Theta_1, \Theta_2, \Theta_3, \Theta_4)\in (0, 1/4)$, $\alpha_1=\alpha_1(d, \beta,\Theta_1, \Theta_2, \Theta_3, \Theta_4)\in (0,1)$  and $C_1=C_1(d, \beta,\Theta_1, \Theta_2, \Theta_3, \Theta_4, \Theta_5, m)>0$ such that for any $x_0\in \R^d,r\in (0,\wt{r}_0)$, any bounded function $u$ defined in $\R^d$ that is harmonic with respect to $\mathcal{G}$ in $B(x_0, r)$, 
 \begin{equation*} 
 \|u(x)-u(y)\|\leq C_1\|u\|_{\infty}(\frac{\|x-y\|}{r})^{\alpha_1}, \quad \hbox{ for any } x,y \in B(x_0, r/2).
 \end{equation*}
 Here $\|\cdot \|$ is the supremum norm of a vector.
 \end{theorem}

 The next one is on the scale-invariant Harnack inequality for nonnegative $\mathcal{G}$-harmonic functions.
 
   For any nonnegative function $u$ defined in $\R^d$ that is $\mathcal{G}$-harmonic in $B(x_0, r)$, $x_0\in \R^d, r>0$, define
\begin{equation}\label{defofh}
h(x,k):=\bE^{(x,k)}[u(X_{\tau_{B(x_0, r/2)}}, \Lambda_{\tau_{B(x_0, r/2)}});\tau_{B(x_0,r/2)}<\tau_1], \quad (x,k)\in B(x_0,r/2),
\end{equation}
where
\begin{equation}\label{defoftau1}
 \tau_1=\inf\{t>0:\Lambda_t\neq \Lambda_0\},
 \end{equation}
 is the first switching time for $(X,\Lambda).$

  \begin{theorem}\label{T:1} Let $\mathcal{G}\in \sN(c_0, m,\gamma, \beta, \Theta_1, \Theta_2, \Theta_3, \Theta_4,\Theta_5, c_1, \vartheta)$ with $Q$ satisfying \eqref{condQ}.Then there exist constants $0\leq \wt{r}_1<1$ and $C_2,C_3>0$, which  depends only on \hfill\break\ $(d,m, c_0,\beta, \gamma, \Theta_1,\Theta_2,\Theta_3,\Theta_4, \vartheta, Q^0)$, such that for any $x_0
 \in \R^d, r\in (0,\wt{r}_1]$, any nonnegative function $u$ defined on $\R^d\times \sM$ that is  $\mathcal{G}$-harmonic on $B(x_0, r)\times \sM$, we have for any $x,y \in B(x_0, r/8) ,k\in \sM,$ 
 
$$C_2\Big(h(y,k)+\sum_{l\in E(k)}r^{m_{kl}}h(y,l)\Big)\leq u(x,k)\leq C_3 \Big(h(y,k)+\sum_{l\in E(k)}r^{m_{kl}}h(y,l)\Big),$$
where $h(y,l)$ is defined in \eqref{defofh},
$E(k)=\{l\in \sM\setminus \{k\}:\Psi(n,k,l)\neq \emptyset,\hbox{ for some }n\in \N\}$ and 
$m_{kl}$ is the smallest integer $n$ such that $\Psi^0(n,k,l)\neq \emptyset$.

In particular, for each $k\in \sM$
$$ u(x,k)\leq Cu(y,k),\hbox { for any }x,y \in B(x_0, r/8),$$
where $C=C_3/C_2$.
\end{theorem}
Under the assumption of strict irreducibility for $\mathcal{G}$, we have the following scale-invariant full rank Harnack inequality for nonnegative $\mathcal{G}$-harmonic functions.
\begin{theorem}\label{T:2}(Full Rank Harnack Inequality) Let $\mathcal{G}\in\sN(c_0, m,\gamma, \beta, \Theta_1, \Theta_2, \Theta_3, \Theta_4,\Theta_5, c_1, \vartheta)$ with $Q$ satisfying \eqref{condQ} and $\mathcal{G}$ is strictly irreducible. Then there exist constants $\wt{C}_1,\wt{C}_3>0$ and  $0\leq \wt{r}_1<1$  depending on $(d, \beta, \gamma, \Theta_1,\Theta_2,\Theta_3,\Theta_4, \vartheta,m,c_0, q_0,Q^0)$,
such that for any $r\in (0,\wt{r}_1]$, any nonnegative function $u$ defined on $\R^d\times \sM$ that is  $\mathcal{G}$-harmonic on $B(x_0,r)\times \sM$, we have
$$u(x,k)\leq \wt{C}_1(u(y,l)+h(y,k))\quad \hbox{ for any }(x,k),(y,l)\hbox{ in } B(x_0,r/8)\times \sM.$$
 In particular,
$$u(x,k)\leq \wt{C}_3r^{-2} u(y,l) \quad \hbox{ for any }(x,k),(y,l)\hbox{ in } B(x_0,r/8)\times \sM \hbox{ with }k\neq j. $$
Here  $h(y,k)$is defined in \eqref{defofh},  and $\wt{r}_1, q_0$ are the constants in Theorem \ref{T:1} and \eqref{lowerbdofQ0} respecitvely.
\end{theorem}
 The rest of the paper is organized as follows. In Section 2, we derive the two-sided scale invariant Green function estimate and the Martin integral representation formula for nonnegative harmonic functions with respect to the operator $\sL_i+q_{ii}$, as well as the probabilistic representation formula for nonnegative $\mathcal{G}$-harmonic functions in term of harmonic functions for $\sL_i+q_{ii}, i\in \sM$ in small balls.  The proof of the H\"older regularity of bounded $\mathcal{G}$-harmonic functions is given in Section 3. The proof of Theorem \ref{T:1} is given in Section 4 and the proof for Theorem \ref{T:2} is shown in Section 5.

 \section {Preliminaries}\label{S:1} 
A Borel measurable function $f$ is said to be in {\bf{Kato class}} $\mathbb{K}_d$ if and only if for each ball $B$ in $\R^d$:
\begin{equation}\label{katoclass}
\lim_{r\rightarrow 0}[\sup_{x\in B}\int_{|y-x|\leq r}|f(y)|g(y-x)dy]=0,
\end{equation}
where
\begin{equation}\label{functg}
g(x)=
\begin{cases}
-\ln|x| \hbox{ when } d=2;\\
|x|^{2-d} \hbox{ when } d\geq 3.\\
\end{cases}
\end{equation}
\begin{remark}
The boundedness assumption of $q_{kk}$ in \eqref{conditiononQ} implies that each $q_{kk}\in \mathbb{K}_d, k\in\sM$.
\end{remark}
 
The following rough scaling property of the infinitesimal generator identifies the class of the operators for which the scale-invariant Harnack inequalities  will be established in this paper. 
\begin{lemma}\label{scaling}
Suppose that $\{(X_t,\Lambda_t), t\geq 0; \bP^{(x,i)}, (x,i)\in \R^d\times \sM\}$ is a Hunt process having the infinitesimal generator  $\mathcal{G} \in \sN(c_0, m,\gamma, \beta, \Theta_1, \Theta_2, \Theta_3, \Theta_4,\Theta_5, c_1, \vartheta)$. 
For any $\lambda \in (0, 1]$, define
$Y_t=\lambda^{-1} X_{\lambda^2 t}$, $\Gamma_t=\Lambda_{\lambda^2 t}$, and $\fP^{(x,i)}= \bP^{(\lambda x, i)}$. Then
$\{(Y_t, \Gamma_t), t\geq 0; \fP^{(x,i)}, (x,i)\in \R^d\times \sM\}$ is a Hunt process corresponding to the infinitesimal generator
\begin{eqnarray*}
\mathcal{G}^{(\lambda)} f(x,i) =\sL_{i}^{(\lambda)}f(x,i)+Q^{(\lambda)}(x)f(x,\cdot)(i),\quad \hbox{ for }(x,i)\in \R^d\times \sM,
\end{eqnarray*}
where $f(\cdot, i)\in C^2_b(\R^d),i\in\sM$,
\begin{eqnarray}\label{scaledoperator}
\sL_i^{(\lambda)} f(x,i) &=& \frac{1}{2}\sum_{k,l=1}^{d} a_{k,l}(\lambda x,i)\frac{\partial ^2 f }{\partial x_k \partial x_l} (x,i)\nonumber\\
&&+\Big(\lambda b_1(\lambda x,i) -\int_{\R^d\setminus\{0\}}z \1_{\{1< |z|\leq 1/\lambda\}} \lambda^{d+2}b_2(\lambda x, \lambda z,i)j_0(\lambda z)dz\Big)\cdot \nabla f(x,i) \nonumber\\
&&  +\int_{\R^d\setminus\{0\}} \left( f(x+z,i)-f(x,i)-  z \1_{\{|z|\leq 1\}}  \cdot \nabla f(x,i) \right) \lambda^{d+2} b_2(\lambda x, \lambda z,i)j_i(\lambda z) dz\nonumber\\
&& +\lambda^2Q(\lambda x)f(x,\cdot)(i),
\end{eqnarray}
and
$$q^{(\lambda)}_{ij}(x)=\lambda^2q_{ij}(\lambda x), \hbox{ for every }i,j\in \sM.$$
In particular, $\mathcal{G}^{(\lambda)} \in \mathcal{N}(\Theta_1,\Theta_2+\Theta_3\Theta_4,  \Theta_3, \Theta_4,\Theta_5,\vartheta, \beta, \gamma ,c_1)$
with $Q^{0,(\lambda)}=\lambda^2 Q^0 $ in \eqref{condQ} for every $\lambda \in (0, 1]$. 
\end{lemma}

\proof  For any  $f\in C_b^2(\R^d, \sM)$, define $f_{\lambda}(x,i)=f(x/\lambda,i)$. Denote by $\fP^{(x,i)}$ the transition semigroup of  $(Y,\Gamma)$.Then
$$\fP_tf(x,i)=\fE^{(x,i)}[f(Y_t, \Gamma_t)]=\bE^{(\lambda x,i)}f_{\lambda}[(X_{\lambda^2t}, \Lambda_{\lambda^2t})]=\bP^{(\lambda x,i)}_{\lambda^2t}f_{\lambda}(\lambda x, i).$$
Then
\begin{eqnarray}
&&\mathcal{G}^{(\lambda)}f(x,i)\nonumber\\
&=&\lim_{t\rightarrow 0}\frac{\fE^{(x,i)}[f(Y_t, \Gamma_t)]-f(x,i)]}{t}\nonumber\\
&=&\lim_{t\rightarrow 0}\frac{\bE^{(\lambda x,i)}[f_{\lambda}(X_{\lambda^2t}, \Lambda_{\lambda^2 t})]-f_{\lambda}(x,i)]}{\lambda^2t}\lambda^2=\lambda^2\mathcal{G}f_{\lambda}(\lambda x,i).\nonumber\\
&=& \frac{1}{2}\sum_{i,j=1}^{d} a_{i,j}(\lambda x,i)\frac{\partial ^2 f }{\partial x_i \partial x_j} (x,i)
+\Big(\lambda b_1(\lambda x,i) -\int_{\R^d\setminus\{0\}} z \1_{\{1< |z|\leq 1/\lambda\}} \lambda^{d+2} b_2(\lambda x, \lambda z,i)j_0(\lambda z)dz)\Big)\cdot \nabla f(x,i) \nonumber\\
&&  +\int_{\R^d\setminus\{0\}} \left( f(x+z,i)-f(x,i)-z \1_{\{|z|\leq 1\}}  \cdot \nabla f(x,i) \right)  \lambda^{d+2} b_2(\lambda x, \lambda z,i)j_i(\lambda z) dz +\lambda^2Q(\lambda x)f(x,\cdot)(i)\nonumber\\
&=& \sL^{(\lambda)} f(x,i)+Q^{(\lambda)}f(x,i).
\end{eqnarray}
Hence the process $(Y,\Gamma)$ has diffusion matrix $a(\lambda x)$ which satisfies the same uniform ellipticity and H\"older continuity condition
\eqref{e:2.4} as $a(x)$, the drift coefficient vector
\begin{eqnarray}
|b_1^\lambda (x,i)|:&=& |\lambda b_1(\lambda x,i)-\lambda \int_{\R^d\setminus\{0\}} w \1_{\{\lambda \leq |w|\leq 1\}}  b_2(\lambda x, w,i)j_0(w)d(w)|\nonumber\\
&\leq& \Theta_2+\Theta_3 \int_{\R^d\setminus\{0\}}|w|^2 \1_{\{\lambda \leq |w|\leq 1\}}\frac{c_1}{|w|^{d+\beta}}dw\nonumber\\
&\leq& \Theta_2+\Theta_3\Theta_4
\end{eqnarray}
and by\eqref{UJScondition}, the jumping kernel $J^{b_2}_\lambda (x, x+z):=  \lambda^{2} b_2(\lambda x, \lambda z)j_0(\lambda z)$ satisfies
\begin{eqnarray}
J^{b_2}_\lambda (x, x+z)&=& \lambda^{2} b_2(\lambda x, \lambda z)j_0(\lambda z)\nonumber\\
&\leq& \frac{\vartheta}{|B(\lambda x,\lambda r)|}\int_{B(\lambda x,\lambda r)}\lambda^{2} b_2(\lambda u, \lambda (x+z))j_0(\lambda z)d(\lambda u)\nonumber\\
&=&\frac{\lambda^d\vartheta}{|B(\lambda x,\lambda r)|}\int_{B(x,r)}J^{b_2}_{\lambda} (u, x+z)du\nonumber\\
&=&\frac{\vartheta}{|B(x, r)|}\int_{B(x, r)}J^{b_2}_{\lambda} (u, x+z)du,
\end{eqnarray}
which satisfies the UJS condition with the same constant $\vartheta$.
Clearly for $\lambda \in (0, 1]$, 

 $$
\|b_2^{\lambda}\|_{\infty}:= \lambda^{2}\| b_2(\lambda x, \lambda z,i)\|_{\infty}\leq \Theta_3,
 $$
 and 
 $$ 
 j_i^{\lambda}(z):= j_i(\lambda z)\leq \frac{c_1}{|\lambda z|^{d+\beta}},\quad \hbox{ for } |z|\leq 1/\lambda.
$$
with $$
\int_{\R^d\setminus \{0\}}(1\wedge |z|^2)j_i^{\lambda}(z)d z=\int_{\R^d\setminus 0}(\lambda^2 \wedge |w|^2)j_i(w)dw\leq \Theta_4.
$$
 Also, by \eqref{conditiononQ} and \eqref{lowerbdofQ0}, since 
\begin{equation*}
 c_0\lambda^2q_0\leq c_0\lambda^2q_{ij}^0\leq q^{(\lambda)}_{ij}(x)\le\lambda^2 q_{ij}^0\leq  \lambda^2\Theta_5\leq\Theta_5  \quad a.e. \hbox{ in }\R^d.
\end{equation*}
Then $Q^{0,(\lambda)}=\lambda^2 Q^0$ with $c_0^{(\lambda)}=c_0$.
This shows that $\mathcal{G}^{(\lambda)} \in \mathcal{N}(\Theta_1,\Theta_2+\Theta_3\Theta_4, \Theta_3, \Theta_4,\Theta_5, \vartheta, \beta, \gamma)$  with 
$Q^{0,(\lambda)}=\lambda^2 Q^0$ for every $\lambda \in (0, 1]$.\\

\qed

Suppose that $\{(X_t,\Lambda_t), t\geq 0; \bP^{(x_0,i)}, (x_0,i)\in \R^d\times \sM\}$ is a Hunt process having the infinitesimal generator  $\mathcal{G} \in \mathcal{N}(\Theta_1,\Theta_2,  \Theta_3, \Theta_4,  \Theta_5, \gamma, \beta, c_1)$ starting from $(x_0, i)$, define
 \begin{equation*}\label{subprocessofX} 
  X^0_t=
 \begin{cases}
X_t, \quad  t<\tau_1;\\
\partial, \quad  t\geq \tau_1;\\
 \end{cases}
  \end{equation*}
where  $\tau_1$ is defined in \eqref{defoftau1} and $\partial$ is a cemetery point. Then by Lemma 3.7 in \cite{wang},
\begin{eqnarray}\label{MP1}
&&\{X_t^{0}, t\leq \tau_1; \bP^{(x_0,i)},x\in \R^d\} \text{  is a Hunt process corresponding to the operator }\nonumber\\
&&\sL_i+q_{ii} \hbox{ starting from } x. 
\end{eqnarray}

 For any $\varphi\geq 0, \alpha\geq 0$, define 
  \begin{eqnarray}\label{alpharesolvent}
  \wt{G}^i_\alpha\varphi(x):&=&\bE^{(x,i)}[\int_0^\infty e^{-\alpha s}\varphi(X^{0}_s)ds]\nonumber\\
  &=&\bE^{(x,i)}[\int_0^{\tau_1} e^{-\alpha s}\varphi(X_s)ds],
    \end{eqnarray}

    \begin{proposition}\label{distribution1}
    Suppose that $\{(X_t,\Lambda_t), t\geq 0; \bP^{(x,i)}, x\in \R^d,i \in \sM\}$ is the 
 Hunt  process having the infinitesimal generator  $\mathcal{G} \in \sN(c_0, m,\gamma, \beta, \Theta_1, \Theta_2, \Theta_3, \Theta_4,\Theta_5, c_1, \vartheta)$ starting from $(x, i)$. Then for every nonnegative Borel measurable function  $\varphi$ defined on $\R^d$, $x \in \R^d$, any $\alpha\geq 0$, 
   $$ \bE^{(x,i)}[e^{-\alpha \tau_1}\varphi(X_{{\tau_1}^{-}})]=\wt{G}_\alpha^i(-q_{ii}\varphi)(x),  $$
   where $\tau_1$ is defined in \eqref{defoftau1}, and $\wt{G}_\alpha^i$ is defined in \eqref{alpharesolvent}.
    \end{proposition}
    \proof 
Let $\{Z_t^{i},t\leq \xi;\fP^x,x\in \R^d\}$ be the Hunt process corresponding to the operator $(\sL_i, C^2_b(\R^d))$. Define 
 \begin{equation}\label{defofxi}
 \xi=\inf\{t>0:-\int_0^tq_{ii}(Z_s^i)ds>\eta\},
 \end{equation}
  where $\eta$ is an exponential random variable with mean $1. $
Let the  probability measure $\bP^x$ on $\Omega$ determined by 
\begin{equation}\label{defintionofP}
\bE^{x}[f(Z_t^i)]=\fE^x[f(Z_t^{i});t<\xi]=\fE^x[e_q(t)f(Z_t^{i})],\quad t\geq 0,
\end{equation}
for every $f\in C^2_b(\R^d)$, where $e_q(t):=\exp(\int_0^t q_{ii}(Z_s^{i})ds).$
Then by Ito's formula,
\begin{eqnarray}\label{MP2}
&& \{Z^i_t, t\leq \xi; \bP^{x}, x\in \R^d \}\text{ is a Hunt process corresponding to the operator }\nonumber\\
&&\sL_i+q_{ii} \hbox{ starting from }x, \hbox{ whose lifetime is }\xi.
\end{eqnarray}
Then by Theorem 5.2 \cite{komatsu},  \eqref{MP1}, \eqref{MP2}, $\bE^{(x,i)}[\mathds{1}_{\{0<t\leq \tau_1\}}f(X^0_t)g(\tau_1)]=\bE^x[\mathds{1}_{\{0<t\leq \xi\}}f(Z^i_t)g(\xi)]$ for any $f,g\in C^2_b(\R^d)$.Then  
 for any $\varphi\geq 0, x\in \R^d$,
\begin{equation*}
\bE^{(x,i)}[e^{-\alpha\tau_1}\varphi(X_{\tau_1^{-}})]=\bE^{(x,i)}[e^{-\alpha\tau_1}\varphi(X^{0}_{\tau_1^{-}})]=\bE^x[e^{-\alpha\xi}\varphi(Z^{i}_{\xi^{-}})]
\end{equation*}
and 
\begin{eqnarray}\label{alpharesolventequi}
\wt{G}_\alpha^i\varphi(x)=\bE^{(x,i)}[\int_{0}^\infty e^{-\alpha t}\varphi(X^0_t)dt]=\bE^{x}[\int_{0}^\infty e^{-\alpha t}\varphi(Z^i_t)dt].
\end{eqnarray}
Let $\mathcal{F}_t=\sigma\{Z^i_s;0\leq s\leq t\}$.Then by \eqref{defofxi},\eqref{defintionofP},\eqref {alpharesolventequi}, and  (61.2) on \cite[p.286]{General} (by putting $m_t=e^{\int_0^tq_{ii}(Z^i_s)ds}\mathds{1}_{\{t<\xi\}}$ there and noticing that $\fP^{x}-a.e. Z_{\xi^{-}}=\partial$) that 
\begin{eqnarray*}
\bE^x[e^{-\alpha\xi}\varphi(Z^{i}_{\xi^{-}})]&=&\fE^x[\int_0^\infty e^{-\alpha t}\varphi(Z^{i}_{t})(-dm_t)]\nonumber\\
&=&\fE^{x}[\int_0^\infty e^{-\alpha t}\varphi(Z^i_t) (-q_{ii}(Z^{i}_t))\exp(\int_{0}^t q_{ii}(Z^{i}_s) ds)) dt]\nonumber\\
&=&\bE^{x}[\int_0^\infty e^{-\alpha t}(-q_{ii}\varphi)(Z^{i}_t)dt]\nonumber\\
&=&\wt{G}^i_\alpha(-q_{ii}\varphi)(x).
\end{eqnarray*}
\qed
 \begin{remark}
By \eqref{alpharesolventequi}, when $\alpha=0$, we know that $\wt{G}^{i}$ is the Green operator of $\sL_i+q_{ii}$  in $\R^d$ with zero Dirichlet boundary condition.
\end{remark}
For any nonnegative Borel measurable  function $\phi$ defined on $\R^d\times \sM$, define
  \begin{equation}\label{defintionofu}
  u(x,i):=\bE^{(x,i)}[\phi(X_{\tau_D},\Lambda_{\tau_D})] \hbox{ for any }(x,i)\in D\times \sM.
 \end{equation} 

\qed

    \begin{proposition}\label{representofu}   Let $D$ be a bounded and connected open set in $\R^d$.
    For any nonnegative Borel measurable  function $\phi$ defined on $\R^d\times \sM$,  the function $u$ defined in \eqref{defintionofu} satisfies
\begin{equation}\label{represenofu1}
u(x,i)=h(x,i)+\sum_{j=1,j\neq i}^m\wt{G}_D^i(q_{ij}u(\cdot, j))(x), \hbox{ for every } (x, i)\in D\times \sM,
\end{equation}
where 
\begin{equation}\label{defofhi}
h(x,i)=\bE^{(x,i)}[\phi(X_{\tau_D}, \Lambda_{\tau_D}); \tau_D< \tau_1],
\end{equation}
$\tau_D=\inf\{t>0:(X_t,\Lambda_t)\notin D\times \sM\}$ and
$ \wt{G}^i_D$ is the Green operator with respect to the operator $\mathcal{L}_i+q_{ii}$ on $D$ with zero Dirichlet boundary condition.
  \end{proposition} 
  
  \proof 
  First consider any nonnegative function $\phi\in C_c^{\infty}(\R^d\times \sM)$, where  $C_c^{\infty}(\R^d\times \sM)$ is the space of functions $f$ such that  $f(\cdot,i)\in C_c^{\infty}(\R^d)$ for each $i\in \sM$. Then by \eqref{defofhi},
  \begin{eqnarray}\label{part1}
  u(x,i)&=&\bE^{(x,i)}[\phi(X_{\tau_D}, \Lambda_{\tau_D});\tau_D< \tau_1]+\bE^{(x,i)}[\phi(X_{\tau_D}, \Lambda_{\tau_D});\tau_D> \tau_1]+\bE^{(x,i)}[\phi(X_{\tau_D}, \Lambda_{\tau_D});\tau_D= \tau_1]\nonumber\\
  &=&h(x,i)+I+II;
   \end{eqnarray}
   Then by the strong Markov property of  $(X, \Lambda)$, Lemma 3.6 in \cite{wang}, and Proposition \ref{distribution1}, 
    \begin{eqnarray}\label{part2}
    I&=&\bE^{(x,i)}[\bE^{(X_{\tau_1}, \Lambda_{\tau_1})}[\phi(X_{\tau_D}, \Lambda_{\tau_D})];\tau_D>\tau_1]\nonumber\\
    &=&\bE^{(x,i)}[u(X_{\tau_1}, \Lambda_{\tau_1});\tau_D>\tau_1]\nonumber\\
    &=&\sum_{j=1,j\neq i}^{m}\bE^{(x,i)}[\mathds{1}_{\{\tau_1<\tau_D\}}u(X_{{\tau_1}^{-}},j)(-\frac{q_{ij}}{q_{ii}}(X_{{\tau_1}^{-}}))]\nonumber\\
    &=&\sum_{j=1,j\neq i}^{m}\bE^{(x,i)}[\int_0^{\tau_1} \mathds{1}_{\{\tau_1<\tau_D\}}u(X_s,j)\frac{q_{ij}}{q_{ii}}(X_s) q_{ii}(X_s)ds)]\nonumber\\
   &=& \sum_{j=1,j\neq i}^{m}\bE^{(x,i)}[\int_0^{\tau_D\wedge \tau_1}\mathds{1}_{\{\tau_1<\tau_D\}}u(X_s,j)\frac{q_{ij}}{q_{ii}}(X_s) q_{ii}(X_s)ds)]\nonumber\\
     &=& \sum_{j=1,j\neq i}^{m}\bE^{(x,i)}[\int_0^{\tau_D}u(X_s^0,j)\frac{q_{ij}}{q_{ii}}(X_s^0) q_{ii}(X_s^0)ds)]\nonumber\\
    &=&\sum_{j=1,j\neq i}^{m}\wt{G}^{i}_D(q_{ij}u(\cdot,j))(x).
       \end{eqnarray}
       By the definition of $\tau_1$ in \eqref{defoftau1}, $\tau_1$ is $\mathcal{F}_{\tau_1^{-}}$measurable, then again by Lemma 3.6 in \cite{wang}, 
 and Proposition \ref{distribution1},
           \begin{eqnarray}\label{part3}
           II&=&\bE^{(x,i)}[\mathds{1}_{D^c}(X_{\tau_1})\mathds{1}_{\{\tau_D=\tau_1\}}\phi(X_{\tau_1}, \Lambda_{\tau_1})]\nonumber\\
           &=&\sum_{j=1,j\neq i}^{m}\bE^{(x,i)}[\mathds{1}_{D^c}(X_{{\tau_D}^{-}})\mathds{1}_{\{\tau_D=\tau_1\}}\phi(X_{\tau_1^{-}}, j)\frac{q_{ij}}{q_{ii}}(X_{\tau_1^{-}})]\nonumber\\
           &=&0.
             \end{eqnarray}
             Then by \eqref{part1},\eqref{part2},\eqref{part3},
            \begin{eqnarray*}   
         u(x,i)&=&h(x,i)+\sum_{j=1,j\neq i}^{m} \wt{G}^{i}_D(q_{ij}u(\cdot,j))(x).
         \end{eqnarray*}
         
        Finally, by monotone convergence theorem, \eqref{represenofu1} holds for any nonnegative Borel measurable function $\phi$ defined on $\R^d\times \sM$ harmonic with respect to $\mathcal{G}$ in $D$.
        
      \qed
      
    Recall that an open set $D$ in $\R^d$ (when $d\geq 2$) is said to be $C^{1,1}$ if there exist
a localization radius $R_0>0$ and a constant $\Xi_0>0$ such that for every $Q\in \partial D$ , there exists a $C^{1,1}$ function $\phi=\phi_Q:\R^{d-1}\rightarrow \R$ satisfying $\phi(0)=\nabla \phi(0)=0$, $\|\phi\|_{\infty}\leq \Xi_0,|\nabla \phi(x)-\nabla \phi(y)|\leq \Xi_0|x-y|$, and an orthonormal coordinate system $CS_Q:y=(y_1,...y_{d-1},y_d)=:(\wt{y},y_d)\in \R^{d-1}\times \R$ with its origin at $Q$ such that
$$B(Q,R_0)\cap D=\{y=(\wt{y},y_d)\in B(0, R_0) \hbox{ in } CS_Q: y_d>\phi(\wt{y})\}.$$
The pair $(R_0,\Xi_0)$ is called the characteristics of the $C^{1,1}$ open set $D$.
   
 For any $\sL\in \sN(\Theta_1, \Theta_2, \Theta_3,\Theta_4,\gamma, \beta, c_1)$, let  $(X_t,t\geq 0;\bP^x, x\in \R^d)$ be the process having the infinitesimal generator $\sL$ and denote $X^D$ as the subprocess of $X$ in $ D$. That is, $X_t^D(\omega)=X_t(\omega)$ if $t<\tau_D(\omega)$ and $X^D(\omega)=\partial$ if $t\geq \tau_D(\omega)$, where $\partial$ is a cemtery state. If $p(t,x,y)$ is the transition density of the process $X$. Then 
 $$p_D(t,x,y)=p(t,x,y)-\bE^x[p(t-\tau_D,X_{\tau_D},y),\tau_D<t]$$ is the transition density of $X^D$. Let {\bf{Green function}} of $\sL$ in a bounded $C^{1,1}$ domain $D$. That is, $G_D(x,y)=\int_{0}^{\tau_D}p_D(t,x,y)dt $ for $x,y\in D$. 
Then for any $\lambda \in (0, 1)$, we denote the Green function for $ \sL$ in $ \lambda D$  as  $G_{ \lambda D}(x, y)$, for $x,y\in \lambda D$ and the the Green function for $\sL^{(\lambda)}$ of the form \eqref{scaledoperator} in $D$ as $G_{D}^{(\lambda)}(x,y)$, for $x, y \in D$.

Let $G^{\Delta}_{\lambda D}$ be the Green function for  $\Delta$ in  $\lambda D$.
Notice that by Brownian scaling, $$G^{\Delta}_{\lambda D}(x,y)=\lambda^{2-d}G^{\Delta}_D(x/\lambda, y/\lambda), \hbox{ for every }x,y\in \lambda D,x\ne y.$$

For the details of how to transform the relationship between $G_{D}^{(\lambda)}(x/\lambda,y/\lambda)$ and $G^{\Delta}_D(x/\lambda, y/\lambda)$ to the relationship between  $G_{\lambda D}(x, y)$ and $G^{\Delta}_{\lambda D}(x, y), x, y
\in \lambda D$, the reader can refer to \cite{BHP}. Next we quote a series of theorems for the sharpe two-sided scale-invariant Green function estimates from \cite{BHP}, which helps us obtain scale-invariant properties for the Green function of the killed operator $\sL+q$ as a preparation to derive the Martin integral representation formula for nonnegative harmonic functions with respect to  $\sL+q$, where $q<0$ and $\|q\|_{\infty}=\Theta_5$.

The following lemmas on the scale-invariant comparison inequalities of the Green functions for the Laplacian operator $\Delta$ in small balls  can be founded in \cite[(2.19)]{BHP} for $d\geq 3$ and in \cite[Lemma 3]{greenfuncd2} for $d=2$ respectively.

 \begin{lemma}\label{2.19}
 There exists a constant $K_1=K_1(d)>1$ for $d\geq 3$ such that for any $x_0\in \R^d,r>0$, such that for any $x\neq y\in B(x_0,r), x\neq y$,
  \begin{eqnarray}\label{greenfunc1}
  &&\frac{K_1^{-1}}{|x-y|^{d-2}}\Big(1\wedge \frac{\delta_{B(x_0, r)}(x)}{|x-y|}\Big)\Big(1\wedge \frac{\delta_{B(x_0, r)}(y)}{|x-y|}\Big)\leq G_{B(x_0, r)}^{\Delta}(x,y)\nonumber\\
  &\leq& \frac{K_1}{|x-y|^{d-2}}\Big(1\wedge \frac{\delta_{B(x_0, r)}(x)}{|x-y|}\Big)\Big(1\wedge \frac{\delta_{B(x_0, r)}(y)}{|x-y|}\Big),
  \end{eqnarray}
  where $\delta_{B(x_0, r)}(y)$ is the Euclidean distance between the point $y$ to $B(x_0, y)^c$.
\end{lemma}
\begin{lemma}\label{2.9}
When $d=2$, for any $x_0\in \R^d, r\in (0,1)$, any $x,y\in B(x_0, r)$, we have
\begin{equation}\label{greenfunc2}
\frac{1}{2\pi}\ln\Big(1+\frac{\delta_{B(x_0, r)}(x)\delta_{B(x_0, r)}(y)}{|x-y|^2})\leq G^{\Delta}_{B(x_0, r)}(x,y)\Big)\leq\frac{1}{2\pi}\ln\Big(1+4\frac{\delta_{B(x_0, r)}(x)\delta_{B(x_0, r)}(y)}{|x-y|^2}\Big),
\end{equation}
\end{lemma}
We also quote the two-sided scale-invariant Green function estimate for each $\sL$  from \cite[Lemma 2.17]{BHP}.
\begin{lemma}\label{lemma217}
Let $\sL \in \sN (\Theta_1, \Theta_2, \Theta_3,\Theta_4,\gamma, \beta, c_1)$ and $D$ be a bounded $C^{1,1}$ domain with characteristics $(R_0, \Xi_0)$.There exist positive constants $\delta_1\in (0,1)$depending on \hfill\break $(d, \beta, \gamma, \Theta_1, \Theta_2,\Theta_3, \Theta_4,R_0, \Xi_0,diam(D))$  and  $K_0=K_0(d,\Theta_1,\gamma)>1$ such that for any $\lambda\in (0, \delta_1)$ and bounded functions $b_2$ with $\|b_2\|_{\infty}\leq \Theta_3\vee (\Theta_2+\Theta_3\Theta_4)$,
\begin{equation}\label{greenfunctionbd}
K_0^{-1}G_{\lambda D}^{\Delta}(x,y)\leq  G_{\lambda D}(x,y)\leq K_0G_{ \lambda D}^{\Delta}(x,y) \quad \hbox{ for }x,y\in \lambda D,
\end{equation}
\end{lemma}
\qed

Let  $\sL \in \sN (\Theta_1, \Theta_2, \Theta_3,\Theta_4,\gamma, \beta, c_1)$. Fix $x'\in D$ and let  
\begin{equation}\label{martinkernel}
 M_{D}(x,z)=\frac{G_D(x, z)}{G_D(x', z)}, \hbox{ for }x\in D,y \in D\setminus\{x,x'\}
\end{equation}
 be the {\bf Martin kernel} of $\sL$ for $x\in D, z \in D\setminus\{x,x'\}$. The {\it Martin boundary} for $\sL$ in $D$ is defined to be the set $\partial_MD=D^{\ast}\setminus D$, where $D^{\ast}$ is the smallest compact set for which $M_D(x,z)$ is continuous in $z$ in the extended sense(See 
\cite[Theorem 7.1]{BassPTA} for detailed explanation for this definition). The {\it minimal boundary } $\partial^{min}_MD$ for $\sL$ in $D$ is the collection of all points $z\in \partial_M D$ so that $x\mapsto M_D(x,y)$ is a positive {\it minimal }$\sL$-{\it harmonic function} in $D$ in the sense that if $h\geq 0$ is $\sL$-harmonic in $D$ and $h(x)\leq M_D(x,y)$ on $D$, then $h(x)=cM_D(x,y)$
 for some constant $0<c\leq 1$.
 
By Theorem 3.6 in \cite{BHP} which is quoted below, we can identify the Euclidean boundary as its Martin boundary with respect to the operator $\sL$ of the form \eqref{e:1.3} when $D$ is a bounded $C^{1,1}$ domain.
 
 \begin{theorem}\label{martinbdry}
 Suppose that $D$ is a bounded $C^{1,1}$ domain with characteristics $(R_0, \Xi_0)$.  Then there exists a  positive constant $\delta_1>0$ such that that for any $\lambda\in (0,\delta_1)$, the Martin boundary and minimal Martin boundary $\lambda D$ with respect to $\sL$ of the form \eqref{e:1.3} can be both identified with its Euclidean boundary $\partial (\lambda D)$, where $\delta_1$ is the constant in Lemma \ref{lemma217}.
 \end{theorem}

 Moreover, the Martin integral representation formula for nonnegative harmonic function with respect to $\sL$ is obtained in \cite[Theorem 3.8]{BHP}.
\begin{theorem}
 Let $D$ be a bounded $C^{1,1}$ domain and $\lambda<\delta_1$. If $h$ is a nonnegative harmonic function in $\lambda D$ with respect to $\sL\in \sN(\Theta_1, \Theta_2, \Theta_3,\Theta_4,\gamma, \beta, c_1)$, then there exists a unique measure $\mu_h$ on $\partial (\lambda D)$ such that 

\begin{equation}
h(x)=\int_{\partial (\lambda D)}M_{ \lambda D}(x,z)h(z)\mu_h(dz)+\int_{\overline{ \lambda D}^c}\int_{ \lambda D}G_{ \lambda D}(x,y)J^{b_2}(y,z)dyh(z)dz,
\end{equation}
where  $\delta_1$ is the constant in Lemma \ref{lemma217},  and  $J^{b_2}$ is the jump kernel defined in \eqref{UJScondition}.
\end{theorem}

 Suppose that $\sL \in \sN (\Theta_1, \Theta_2, \Theta_3,\Theta_4,\gamma, \beta, c_1)$ and $q<0$ such that $\|q\|_{\infty}=\Theta_5.$ Let $\wt{G}_{B(x_0,r)}(\cdot, \cdot)$ be the {\bf Green function}  for the operator $\sL+q$ in the ball $B(x_0, r)$.
We are going to prove the existence and a comparison inequality for $\wt{G}_{B(x_0,r)}(\cdot, \cdot)$. First, we quote the $3G$-lemma  from \cite[Lemma 2.1]{weaklycoupled}, which can also be found in \cite{BMgauge, conditionalgauge}.

\begin{lemma}\label{lemma3.14}(3G-Lemma)
For any ball $B\in \R^d$,
\begin{equation}\label{greenfunctioninequality1}
\frac{G^{\Delta}_B(x,y)G^{\Delta}_B(y,z)}{G^{\Delta}_B(x,z)}\leq C_G(g(x-y)+g(y-z))
\end{equation}
 for every $x,y,z\in B$, where the constant $C_G$ depends only on $d$, and the function $g$ is defined in \eqref{functg}.
\end{lemma}
\proof The proof can be found in \cite{BMgauge}.
\qed
\medskip

 For any $\sL \in \sN(\Theta_1, \Theta_2,\Theta_3, \Theta_4, \beta, \gamma,c_1)$, let $\{\bar{X}_t,t\geq 0;\bP^x,x\in\R^d\}$ be the Hunt process corresponding to the operator $\sL$. Define $$\tau_D=\inf\{t\geq 0;\bar{X}_t\notin D\}.$$
 Similar as the $h$-conditioned Brownian motion defined in \cite[Page 131, 132]{conditionalgauge}, for any positive Borel-measurable function $h$, $M_t:=h(\bar{X}_t)/h(\bar{X}_0)$ is a multiplicative functional  $M=\{M_t:0\leq t<\infty \}$ of $\bar{X}$ such that $\bE^{x}[M_t]\leq 1$ for every $t\geq 0$. For any $\mathcal{F}_t$-measurable function $\Phi\geq 0$, where $\mathcal{F}_t=\sigma\{\bar{X}_s:s\leq t\},t>0$, define  
\begin{equation}\label{conditionedexp}
\bE^x_{h}[\Phi;t<\tau_D]:=(h(x))^{-1}\bE^x[\Phi \cdot h(\bar{X}_t);t< \tau_D];\quad x\in D.
\end{equation}
and $\bE^x_{h}[\Phi;t<\tau_D]$ is reduced to  be a probability distribution function in $B$ when $\Phi$ is of the form $\mathds{1}_{\{X_t\in B\}}$, $B\in \mathcal{D},$ where $ \mathcal{D}$ is a Borel $\sigma$-field of $D$.

 When $D$ is a bounded $C^{1,1}$ domain, since the Green function $G_D(\cdot,v)$ for the operator $\sL$ is harmonic with respect to $\sL$ in $D\setminus \{v\}$, then  substituting $h$ with  $G_D(\cdot,v)$ in \eqref{conditionedexp}, we can define the{ \it $G_D(\cdot, v)$-conditioned Markov process} whose state space is $D\setminus\{v\}
\cup\{\partial\}$, and the associated probability and expectation are denoted as $\bP^x_v$ and $\bE^x_v$ respectively. By the definition of Martin kernel in \eqref{martinkernel} and by 
\cite[Lemma 3.2, Lemma 3.4]{BHP}, $M_D(\cdot, z)=\lim_{y\in D, y\rightarrow z}\frac{G_{D}(x, y)}{G_{D}(x', y)}$ exists for any $x \in D, z
\in \partial D$, and it is harmonic with respect to $\sL$ in $D$. Then we can similarly define the {\it $M_D(\cdot, z)$-conditioned Markov process}, whose  associated probability and expectation are denoted as $\bP^x_z$ and $\bE^x_z$ respectively.

\begin{definition}
Let $\sL\in\sN(\Theta_1, \Theta_2,\Theta_3, \Theta_4, \beta,\gamma,c_1)$ and $(\bar{X}_t,t\geq 0;\bP^x,x \in\R^d)$ be the Hunt process corresponding to $\sL$. We say that a set $A\subset \R^d$ is {\bf{nearly Borel}} if for every $x\in \R^d$, there exist Borel subsets $B$ and $B^{'}$ in $\mathcal{B}(\R^d)$ such that $B\subset A\subset B^{'}$ and
$$\bP^x(\bar{X}_t\in B'\setminus B\hbox{ for some }t>0)=0.$$The collection of nearly Borel sets is a $\sigma$-field, which we denote as $\mathcal{B}^n(\R^d)$.We say that a set $A\subset 
\R^d$ is {\bf{polar}} if there exists a set $D\in \mathcal{B}^n(\R^d)$ such that $A\subset D$ and $\bP^x(T_D<\infty)=0$ for every $x\in \R^d$.
\end{definition}
Notice that the definition of nearly Borel set depends on the process $\bar{X}$ and roughly speaking, a set is nearly Borel if the process $\bar{X}$ cannot distinguish it from a Borel set. Also, it is easy to see that $\mathcal{B}(\R^d)\subset \mathcal{B}^n(\R^d)$. We will omit further discussions and the reader can refer to \cite[Page 60]{BG} for more details on nearly Borel sets.
\begin{lemma}\label{polarset}
Let $\sL\in\sN(\Theta_1, \Theta_2,\Theta_3, \Theta_4, \beta,\gamma,c_1)$ and $(\bar{X}_t,t\geq 0;\bP^x,x \in\R^d)$ be the Hunt process corresponding to the operator $\sL$.Then any singleton set is a polar set.
\end{lemma}
\proof 
Fix a point $v\in \R^d$. By Lemma \ref{lemma217} with $D=B(v,1)$, there exists constants $\eps_1=\eps_1(d,\Theta_1, \Theta_2, \Theta_3,\Theta_4, \beta,\gamma )$ and $K_0=(d, \Theta_1,\gamma)$ such that for any $r\in (0,\eps_1)$, \eqref{greenfunctionbd} holds with $\lambda=r$. Since \eqref{greenfunctionbd} is independent of the center point $rv$, then the Green function  $G_{B(v,r)}=G_{r(B(v,1)-v)+v}$ of $\bar{X}$ in $B(v,r)$ satisfies 
\begin{equation}\label{greenfuncbd}
K_0^{-1}G^{\Delta}_{B(v,r)}(x,y)\leq G_{B(v,r)}(x, y)\leq K_0G^{\Delta}_{B(v,r)}(x,y), \hbox{ for }x\in B(v,r), y\in  B(v,r)\setminus\{x\}.
\end{equation}
 Define $\sigma_{B(v, r)}=\inf\{t>0: \bar{X}_t\in B(v, r)\}$ and $\tau_{B(v, r)}=\inf\{t>0: \bar{X}_t\notin B(v, r)\}$. Notice that for any $x \in \R^d$,
 \begin{eqnarray} \label{greeninfty}
U(x,\{v\}):&=&\int_0^\infty \bP_s^{x}(\{v\})ds\nonumber\\
&=&\bE^{x}[\bE^{\bar{X}_{\sigma_{B(v, r)}}}[\int_0^{\tau_{B(v,r)\circ \theta_{\sigma_{B(v, r)}}}}\mathds{1}_{\{v\}}(\bar{X}_s)ds]]+\bE^{x}[\int_{\tau_{B(v, r)}}^{\infty}\mathds{1}_{\{v\}}(\bar{X}_s)ds]\nonumber\\
&=&\bE^x[\int_{\{v\}}G_{B(v,r)}(\bar{X}_{\sigma_{B(v, r)}},y)dy]+\bE^{x}[U(\bar{X}_{\tau_{B(v, r)}},\{v\})].
\end{eqnarray}
Since $\bP_tU(x,\{v\})=\int_0^\infty \bP_{t+s}^{x}(\{v\})ds \rightarrow U(x,\{v\})$ as $t\rightarrow 0$. Then $U(\cdot,\{v\})$ is an excessive function on $\R^d$. 

Denote by $p(t, x,y)$  the transition density function of $\bar{X}_t$. By (2.1) and (2.6) in \cite{BHP}, the transition density of $\bar{X}_t$ for $t>0$, there exists constants $\bar{c}_k,k=1,2,3$ depending on $(d, \beta, \gamma, \Theta_1,\Theta_2,\Theta_3  )$ such that

$$\bar{c}_1^{-1}e^{-c_1t}p_0(t,\bar{c}_2x,\bar{c}_2y)\leq p(t,x, y)\leq \bar{c}_1e^{\bar{c}_1t}p_1(t, \bar{c}_3x,\bar{c}_3y)+t\|J_1\|_{\infty}, $$
where 
$p_0(t, x,y)$ and $p_1(t, x,y)$ are the transition density functions for $\Delta$ and $\Delta^{\beta/2}$ respectively,  and $J_1=\int_{\R^d}\mathds{1}_{\{|z|\geq 1\}}b_2(x,z) j_0(z)dz$.

Then for any $x\in \R^d$,
$\bP_t^{x}(\{v\})=\int_{\{v\}}p(t,x,y)dy=0$ for each $t>0$, and thus $U(x,\{v\})=0$ for every $x\in \R^d$. Also, by \eqref{greenfunc1},\eqref{greenfunc2},\eqref{greenfuncbd}, when $x=v$ in \eqref{greeninfty},
  \begin{equation*}
U(v,\{v\})\geq \int_{\{v\}}G_{B(v,r)}(v,y)dy\geq {K_0}^{-1} \int_{\{v\}}G^{\Delta}_{B(v,r)}(v,y)dy=\infty.
  \end{equation*}
 Then by \cite[Proposition 3.14]{BG},  $\{v\}$ is a polar set.
\qed

Let $\mathcal{F}_{{\tau}^{-}}$ be the $\sigma-$field of events of $Y$ strictly prior to the stopping time $\tau$. That is, $\mathcal{F}_{{\tau}^{-}}$ is $\sigma-$field generated by $\mathcal{F}_0$ and the sets $A\cap \{\tau>t\}$ for $A\in \mathcal{F}_t$ and $t>0$.

 In the following Lemma \ref{estimationforexp}-Theorem \ref{kiillingharmonicrep}, unless we mentioned in particular, we will assume that $\sL\in \sN(\Theta_1, \Theta_2, \Theta_3, \Theta_4, \gamma,\beta,c_1)$ and $(\bar{X}_t,t\geq 0;\bP^{x}, x\in\R^d)$ is the Hunt process 
corresponding to $\sL$  and $\wt{X}$ is the subprocess of $\bar{X}$ with killing rate $-q$, such that $\|q\|_{\infty}=\Theta_5$.
\begin{lemma}\label{estimationforexp}
There exists a constant $r_5=r_5(d,\beta, \gamma, \Theta_1,\Theta_2, \Theta_3, \Theta_4,\Theta_5)<1$ such that for any $r\in (0, r_5), x_0\in \R^d$, we have
\begin{eqnarray}\label{equation5.4}
1/2\leq \bE^{x}_v[\exp(\int_0^{\tau_{B(x_0,r)\setminus \{v\}}}|q|(\bar{X}_s)ds)]\leq 2, \hbox{ for } v\in B(x_0, r),x\in B(x_0, r)\setminus\{v\}, 
\end{eqnarray}
and 
\begin{eqnarray}\label{equation5.10}
1/2\leq \bE^{x}_z[\exp(\int_0^{\tau_{{B}}}|q|(\bar{X}_s)ds)]\leq 2, \hbox{ for } x\in B(x_0,r),z\in \partial B(x_0,r).  
\end{eqnarray}
where   $\tau_{B(x_0, r)\setminus\{v\}}=\inf\{t>0:\bar{X}_t\notin B(x_0, r)\setminus \{v\}\}$, $\tau_B=\inf\{t>0:\bar{X}_t\notin B(x_0, r)\}$, and $\bP^x_v, \bP^x_z$ are defined in \eqref{conditionedexp}.
\end{lemma} 
\proof 
 Fix $x_0\in \R^d$. By Lemma \ref{lemma217} with $D=B(x_0, 1)$, then there exists a constant $\eps_2=\eps_2(d,\beta, \gamma,\Theta_1, \Theta_2, \Theta_3, \Theta_4)>0$ and $K_0= K_0(d,\gamma,\Theta_1)>0$ such that for any $r\in (0,\eps_2 )$, satisfies \eqref{greenfunctionbd} with $\lambda=r$. Since $B(x_0,r)=r(D-x_0)+x_0$, then by denoting $B=B(x_0,r)$, the Green function $G_B(\cdot, \cdot)$ of $\bar{X}$ in $B$ satisfies 
$$K_0^{-1}G^{\Delta}_{B}(x,y)\leq G_{B}(x, y)\leq K_0G^{\Delta}_{B}(x,y), \hbox{ for }x\in B(x_0,r), y\in  B(x_0,r)\setminus\{x\}.$$

Then for any $x\in B,$ $v\in B\setminus \{x\}$,   by Lemma \ref{polarset}, and the definition of the polar set and  the fact that $\{v\}
\in \mathcal{B}^n(\R^d)$,
 $\tau_{B\setminus\{v\}}=\tau_B\quad a.s.$ 
 
Then by the definition of $\bE^x_v$, \eqref{greenfunctionbd},
\begin{eqnarray}\label{equation5.5}
\bE^x_v[\int_0^{\tau_{B\setminus\{v\}}}|q|(\bar{X}_s)ds]&=& (G_B(x,v))^{-1}\bE^x[\int_0^{\tau_B}G_B(\bar{X}_s, v)|q|(\bar{X}_s)ds]\nonumber\\
&=&\int_{B}\frac{G_B(x,y)G_B(y,v)}{G_B(x,v)}|q|(y)dy\nonumber\\
&\leq &K_0^3\int_{B}\frac{G^{\Delta}_B(x,y)G^{\Delta}_B(y,v)}{G^{\Delta}_B(x,v)}|q|(y)dy.
\end{eqnarray}
Then by Lemma \ref{lemma3.14}, there exists a constant $r_5\in (0, \eps_2)$ depending on $(d,K_0,\Theta_5)$ such that for any $r\in (0, r_5)$,
\begin{eqnarray}\label{equation5.6}
K_0^3\int_{B}\frac{G^{\Delta}_B(x,y)G^{\Delta}_B(y,v)}{G^{\Delta}_B(x,v)}|q|(y)dy\leq 1/2,
\end{eqnarray}
 so that by Jenson's inequality, we have for any $ r\in (0,r_5),x\in B(x_0, r), v\in B(x_0, r)\setminus \{x\}$, 
 \begin{equation*}
\bE^x_v[\exp(\int_0^{\tau_{B\setminus\{v\}}}q(\bar{X}_s) ds)]\geq \exp(\bE^x_v[-\int_0^{\tau_{B\setminus\{v\}}}|q|(\bar{X}_s)ds])\geq e^{1/2}> 1/2.
\end{equation*}
Then by \eqref{equation5.5},\eqref{equation5.6} and Khas'minskii's Lemma in \cite[Lemma 3.7]{BMgauge}, we have
\begin{eqnarray*}
\bE^{x}_v[\exp(\int_0^{\tau_{B\setminus\{v\}}}q(\bar{X}_s)ds)]&\leq& \bE^{x}_v[\exp(\int_0^{\tau_{B\setminus\{v\}}}|q|(\bar{X}_s)ds)]\\
&\leq&\frac{1}{1-\sup_{x\in B(x_0,r)}\bE_v^{x}[\int_0^{\tau_{B\setminus\{v\}}}|q|(\bar{X}_s)ds]} \leq 2.
\end{eqnarray*}

 For any $x_0\in \R^d$, by \cite[Lemma 3.2]{BHP} with $D=B(x_0, 1)$ and the fact $B(x_0, r)=r(D-x_0)+x_0$ for any $r\in (0, \eps_2)$, where $\eps_2$ is the constant in the above,  denoting $B=B(x_0,r)$, for any $x\in B,$ $z\in \partial B$,
$M_B(x,z):=\lim_{y\in B,y\rightarrow z}G_B(x,y)/G_B(x',y)$ exists.

Then by the definition of  $\bE^x_z$, and dominated convergence theorem,
\begin{eqnarray}\label{equation5.7}
\bE^x_z[\int_0^{\tau_{ B}}|q|(\bar{X}_s)ds]&=& (M_B(x,y))^{-1}\bE^x[\int_0^{\tau_B}M_B(\bar{X}_s, z)|q|(\bar{X}_s)ds]\nonumber\\
&=&\int_{B}\frac{G_B(x,v)M_B(v,z)}{M_B(x,z)}|q|(v)dv\nonumber\\
&=&\int_{B}\lim_{y\in B, y\rightarrow z}\frac{G_B(x,v)G_B(v,y)}{G_B(x,y)}|q|(v)dv\nonumber\\
&\leq &\lim_{y\in B, y\rightarrow z}\int_{B}\frac{G_B(x,v)G_B(v,y)}{G_B(x,y)}|q|(v)dv.
\end{eqnarray}
Following the above similar argument,  by \eqref{equation5.7}, Lemma \ref{lemma217} with $D=B(x_0,1)$ and $\lambda=r$, Lemma \ref{lemma3.14}, Jenson's inequality and Khas'minskii's Lemma in \cite[Lemma 3.7]{conditionalgauge}, 
for  any $r\in (0, r_5)$,
\begin{eqnarray}\label{equation5.10}
1/2\leq \bE^{x}_z[\exp(\int_0^{\tau_{{B}}}q(\bar{X}_s)ds)]\leq 2, \hbox{ for } x\in B(x_0,r),z\in \partial B(x_0,r).  
\end{eqnarray}
\qed

Next we will establish the lemma on the scale-invariant estimate for the Green function $\wt{G}_{B(x_0, r)}$ of $\sL+q$ in a ball $B(x_0, r)$ for any $ r\in(0,r_5)$ .
\begin{lemma}\label{killgreenfunction} There exists a  constant  $C_0=C_0(d, \gamma,\Theta_1)>1$ such that for any 
$x_0\in \R^d$, $r\in(0, r_5)$, the Green function $\wt{G}_{B(x_0,r)}(\cdot, \cdot)$ for $\sL+q$ in $B(x_0, r)$ exists and satisfies that
\begin{equation}\label{comparision1}
C_0^{-1}G_{B(x_0, r)}^{\Delta}(x,y)\leq \wt{G}_{B(x_0, r)}(x,y)\leq C_0G_{B(x_0, r)}^{\Delta}(x,y), \quad x\in B(x_0, r), y\in B(x_0, r)\setminus\{x\},
\end{equation}
where $r_5$ is the constant in Lemma \ref{estimationforexp} and $C_0=2K_0$ with $K_0$ being the constant in Lemma \ref{lemma217}.
\end{lemma} 
\proof Fix $x_0\in \R^d$.Then by   Lemma 3.5 in \cite{chengauge} with $G(x,y)=G_{B(x_0, r)}(x,y)$ and $\zeta^y=\tau_{B(x_0, r)\setminus{\{y\}}}$,
\begin{eqnarray*}\label{greenfunctionkilled1}
\wt{G}_{B(x_0,r)}(x,y)&=&G_{B(x_0, r)}(x,y) \bE^x_y[\exp(\int_0^{\tau_{B(x_0, r)\setminus{\{y\}}}}|q|(\bar{X}_s)ds)],
\end{eqnarray*}
where $G_{B(x_0,r)}(\cdot,\cdot)$ is the Green function for  $\sL$ in $B(x_0,r)$.
Then the conclusion follows by Lemma \ref{lemma217} and  \eqref{equation5.4}.
\qed
 \begin{corollary}\label{lowerbdexittime}
 For any $\mathcal{G}\in\sN(c_0, m,\gamma, \beta, \Theta_1, \Theta_2, \Theta_3, \Theta_4,\Theta_5, c_1, \vartheta)$, let $((X_t, \Lambda_t),t\geq 0;\bP^{(x,i)},(x,i)\in\R^d\times \sM)$ be the solution solving the martingale problem  $(\mathcal{G}, C_b^2(\R^d\times \sM))$. Then there exists a constant $c_1=c_1(d,\gamma,\Theta_1)>0$ and $c_2=c_2(\eps, d,\gamma,\Theta_1) > 0$ such that 
 for any $x_0\in \R^d, r\in (0,1), \varepsilon \in (0,1),(x,i)\in B(x_0,(1-\varepsilon)r)\times \sM$,
 $$c_2r^2\leq \bE^{(x,i)}[\tau_{B(x_0,r)}\wedge \tau_1]\leq c_1r^2.$$ 
\end{corollary}
\proof Let $\{(X,\Lambda);\bP^{(x,i)}(x,i)\in \R^d\times \sM\}$ be the strong Markov process corresponding to the operator $\mathcal{G}\in\sN(c_0, m,\gamma, \beta, \Theta_1, \Theta_2, \Theta_3, \Theta_4,\Theta_5, c_1, \vartheta)$. And let $\tau^B_{B(x_0, r)}$ be the exit time of the ball $B(x_0, r)$ for a Brownian motion $B$ starting from $x$.   Recall that $\tau_1$ is defined in \eqref{defoftau1}. Then by \eqref{comparision1}, there exist a constant $c_2>0$ depending only on $(\eps, d, \gamma, \Theta_1)$ such that  
\begin{eqnarray}
\bE^{(x,i)}[\tau_{B(x_0, r)}\wedge \tau_1]&=&\int_{\R^d}\wt{G}^{i}_{B(x_0, r)}(x,y)dy\nonumber\\
&\ge& \int_{\R^d}C_0^{-1}G_{B(x_0, r)}^{\Delta}(x,y)dy \nonumber\\
&=&C_0^{-1}\bE^{x}[\tau^B_{B(x_0, r)}]\nonumber\\
&\geq& c_2r^2.
\end{eqnarray}
 
Similarly, by \eqref{comparision1},  there exists a constant $\wt{c}_1>0$ depending only on $( d, \gamma, \Theta_1)$, for any $(x,i)\in B(x_0, r)\times \sM$,
\begin{eqnarray}\label{supexittime}
\bE^{(x,i)}[\tau_{B(x_0, r)}\wedge \tau_1]&\leq& C_0\bE^{x}[\tau^B_{B(x_0, r)}]\nonumber\\
&\leq& \wt{c}_1r^2.
\end{eqnarray}

\qed

Then similar like \cite[Lemma 3.3]{BHP}, by the L\'evy system formula, we have the following lemma.

\begin{lemma}\label{killedrepresenofharmonic}
 For any $x_0\in \R^d, r\in (0,r_5)$, any non-negative function $h$ defined on $\R^d$,  we have for every $x\in B(x_0, r)$,
\begin{equation}
\bE^{x}[h(\wt{X}_{\tau_{B(x_0, r)}}),\wt{X}_{{\tau_{{B(x_0, r)}}}^{-}}\neq \wt{X}_{\tau_{B(x_0, r)}} ]=\int_{\overline{B(x_0, r)}^c}\int_{B(x_0, r)}\wt{G}_{B(x_0, r)}(x,y)J^{b_2}(y,u)dyh(u)du,
\end{equation} 
where $J^{b_2}(y,u)=b_2(y,u-y)j_0(y-u)$ and $r_5$ is the constant in Lemma \ref{killgreenfunction}.
\end{lemma}
 \qed
 
Denote by $M_{B(x_0,r)}^{\Delta}(x,z)$ the Martin kernel of the Laplacian operator $\Delta$ in $B(x_0, r)$, $r\in (0, r_5), x\in B(x_0, r), z\in\partial B(x_0, r)$. In the following, we will obtain the existence and the scale-invariant estimate for the Martin kernel $\wt{M}_{B(x_0,r)}(\cdot, \cdot)$ for $\sL+q$ in $B(x_0, r)$ for any $ r\in (0,r_5)$.

\begin{lemma}\label{killmartinkernel} For any 
$x_0\in \R^d$, $r\in (0,r_5)$, any $x\in B(x_0, r), z\in \partial B(x_0, r)$, $\wt{M}_{B(x_0,r)}(x, z)=\lim_{y\rightarrow z}\frac{\wt{M}_{B(x_0, r)}(x,y)}{\wt{M}_{B(x_0, r)}(x',y)}$ exists and satisfies that
\begin{equation}\label{comparision2}
C_0^{-2}M_{B(x_0, r)}^{\Delta}(x,z)\leq \wt{M}_{B(x_0, r)}(x,z)\leq C_0^2M_{B(x_0, r)}^{\Delta}(x,z), \quad x\in B(x_0, r), z\in \partial B(x_0, r),
\end{equation}
where $r_5$ and $ C_0$ are the constants in Lemma \ref{killgreenfunction}.
\end{lemma} 
\proof The proof of existence of  follows by the definition of $\bE^{x}_y$ in \eqref{conditionedexp} and \cite[Lemma 3.1, Lemma 3.2]{BHP}. And the inequality \eqref{comparision2} follows by Lemma \ref{killgreenfunction}.
\qed

Following the proof of \cite[Lemma 3.4]{CR98}, we can obtain the following lemma on the harmonicity for the killed Martin kernel $\wt{M}_{B(x_0,r)}(\cdot, z)$ with respect to $\sL+q$ in small balls.
\begin{lemma}\label{harmonicmartinkernel}
For any $x_0\in \R^d, r\in (0, r_5)$, for any $z\in \partial B(x_0, r)$,$\wt{M}_{B(x_0, r)}(\cdot,z)$ is harmonic with respect to $\sL+q$ in $B(x_0, r)$.
\end{lemma}

Similar like \cite[Lemma 3.6]{BHP}, we obtain the following theorem.
\begin{theorem}\label{minimalMB}
For any $x_0\in \R^d, r\in (0, r_5)$, the Martin boundary and minimal  Martin boundary of $B(x_0, r)$ for $\sL+q$ can all be identified with its Euclidean  boundary $\partial B(x_0, r)$.
\end{theorem}
\proof The proof follows by Lemma \ref{killmartinkernel} and \cite[Theorem 4.4]{CP02}.\qed

Hence, following the similar arguement in \cite[Theorem 3.8]{BHP}, by Lemma \ref {killgreenfunction}, Lemma \ref{killedrepresenofharmonic} and Theorem \ref{minimalMB}, we have the following theorem.

\begin{theorem}\label{kiillingharmonicrep}
For $x_0\in \R^d, r\in (0,r_5)$, if $\wt{h}$ is a nonnegative function defined in $\R^d$ that is harmonic with with respect to $\sL+q$ in $B(x_0, r)$ , there exists a unique measure $\mu_{\wt{h}}$ on $\partial B(x_0,r)$ such that for any $x\in B(x_0, r),$
\begin{equation}
\wt{h}(x)=\int_{\partial B(x_0,r)}\wt{M}_{B(x_0,r)}(x,z)\wt{h}(z)\mu_{\wt{h}}(dz)+\int_{\overline{B(x_0,r)}^c}\int_{B(x_0,r)}\wt{G}_{B(x_0,r)}(x,y)J^{b_2}(y,z)dy\wt{h}(z)dz.
\end{equation}
\end{theorem} 
Next, given an operator $\mathcal{G}\in \sN (c_0,m,\gamma, \beta,\Theta_1,\Theta_2, \Theta_3, \Theta_4,\Theta_5, c_1,\vartheta)$, for any ball in $\R^d$ with $r(B)\leq r_5$, define $\wt{G}_{B}(\cdot, \cdot)$ on $B\times B$ and another $m\times m$ matrix function $\wt{Q}:\R^d \mapsto \R^{m\times m}$  by 

\begin{eqnarray}\label{GBQ}
 \wt{G}_B(\cdot, \cdot)=
 \left(\begin{matrix}   \wt{G}_B^1(\cdot, \cdot) & 0 & \dots & 0 \\0 & \wt{G}_B^2(\cdot, \cdot)  & \dots & 
0\\ \vdots & \vdots & \ddots & \vdots\\0 & 0 & \cdots & \wt{G}_B^m(\cdot, \cdot) \end{matrix}\right),\quad \quad \wt{Q}=
 \left(\begin{matrix}   0 & q_{12} & \dots & q_{1m} \\ q_{21}  &0  & \dots & 
q_{2m}\\ \vdots & \vdots & \ddots & \vdots\\q_{m1} & q_{m2}& \cdots & 0 \end{matrix}\right),
\end{eqnarray}
where $\wt{G}_B^k(\cdot, \cdot)$ is the Green function of $\sL_k+q_{kk}$ on $B$, and $r_5$ is the constant in Lemma \ref{estimationforexp}. Then $\wt{G}_{B}(\cdot, \cdot)$ is the Green function of the non-coupled operator $\wt{\sL}$,where
\begin{eqnarray}\label{wtL}
 \wt{\sL}=
 \left(\begin{matrix}   \sL_1+q_{11} & 0 & \dots & 0 \\0 &\sL_2+q_{22}  & \dots & 
0\\ \vdots & \vdots & \ddots & \vdots\\0 & 0 & \cdots & \sL_m+q_{mm} \end{matrix}\right).
\end{eqnarray}
 By abusing the notation, let us denote by $\wt{G}_B$ the Green operator of $\wt{\sL}$ and $\wt{G}_B^k$ the Green operator  with respect to $\sL_k+q_{kk}$ on $B$ with zero Dirichlet boundary condition.

To develop a representation of $u$ in terms of $\{h(\cdot,i)\}_{i \in\sM}$ in \eqref{defofhi}, we need the following lemma on the boundedness of the operator $\wt{G}_B\wt{Q}$.
     
 \begin{lemma}\label{norm0}Let $\mathcal{G}\in \sN(c_0, m,\gamma, \beta, \Theta_1, \Theta_2, \Theta_3, \Theta_4,\Theta_5, c_1, \vartheta)$.Then
 there exists a  \hfill \break
 $\delta_3=\delta_3(d, m,\beta, \gamma, \Theta_1,\Theta_2,\Theta_3,\Theta_4,\Theta_5 )>0$ such that for any $x_0\in \R^d, r\in (0, \delta_3]$, 
$\wt{G}_{B(x_0, r)}\wt{Q}$ is a bounded operator on $((B(\overline{B(x_0, r)}\times \sM)),\|\cdot\|_{\infty})$ with 
\begin{equation}
\|\wt{G}_{B(x_0, r)}\wt{Q}\|_{\infty}<1/4.
\end{equation}
 \end{lemma}
 \proof  Fix $x_0\in \R^d$.
   When $d\geq 3$, By Lemma \ref{2.19}, for any $x_0\in \R^d, r>0$,
  \begin{eqnarray} \label{e:3.45}
   \sup_{x\in B(x_0, r)}\int_{B(x_0, r)}G^{\Delta}_{B(x_0, r)}(x,y)dy\leq \omega_dK_1r^2,
         \end{eqnarray}
          where $\omega_d$ is the volume of $d$-dimensional unit sphere.
          
         When $d=2$, since $\ln(1+z)\leq 1+z$ for $z\geq 0$, then by Lemma \ref{2.9}, for any $x_0\in \R^d, r\in (0,1)$,
        \begin{eqnarray} \label{e:3.46}
   \sup_{x\in B(x_0, r)}\int_{B(x_0, r)}G^{\Delta}_{B(x_0, r)}(x,y)dy\leq (\omega_dr^2)(1+4\ln2)/(2\pi)\leq\omega_dK_1r^2 ,
         \end{eqnarray} 
   Let  
   \begin{eqnarray}\label{choiceofr}
   \delta_3=(1/\sqrt{4(m-1)\omega_dC_0(C_G+ K_1)\Theta_5})\wedge (r_5/2),
     \end{eqnarray}
    where $C_G$, $r_5$ and $C_0$ are the constants in Lemma \ref{lemma3.14}, Lemma \ref{estimationforexp} and Lemma \ref{killgreenfunction} respectively.Then for any $r\in (0,\delta_3]$, $u\in (B(B(x_0, r)\times \sM))^m$, it follows that for any $k,l\in \sM$,  by \eqref{represenofu1},\eqref{comparision1}, \eqref{e:3.45},\eqref{e:3.46},
     \begin{eqnarray*}
      \Big|[\wt{G}_{B(x_0, r)}^k\wt{Q}u]_k(x)\Big|&=&\Big|\sum_{l=1,l\neq k}^{l=m}\int_{B(x_0, r)} \wt{G}_{B(x_0,r)}^k(x,y)q_{kl}(y)u(y,l)dy\Big|\\
      &\leq &C_0\Theta_5(m-1)\|u\|_{\infty}\int_{B(x_0,r)}G^\Delta_{B(x_0,r)}(x,y)dy<\|u\|_{\infty}/4.
      \end{eqnarray*}
      Therefore, for any $r\in (0, \delta_3)$, $\wt{G}_{B(x_0, r)}\wt{Q}$ is a bounded operator from $((\mathcal{B}(\overline{B}(x_0, r)\times \sM))^m$ to itself such that 
     $\|\wt{G}_{B(x_0, r)}\wt{Q}\|_{\infty}<1/4$. Here $\mathcal{B}(\overline{B}(x_0, r)\times \sM)$ is the space of bounded functions defined on $\overline{B}(x_0, r)\times \sM$.
     
    \qed  

   \begin{proposition}\label{representofulong} Let $\mathcal{G}\in \sN(c_0, m,\gamma, \beta, \Theta_1, \Theta_2, \Theta_3, \Theta_4,\Theta_5, c_1, \vartheta)$ and $((X_t,\Lambda_t),t\geq 0;\bP^{(x,i)},x\in \R^d\times\sM))$ be the Hunt process corresponding to the operator $\mathcal{G}$.
 For any function $u$ defined in \eqref{defintionofu} with $D=B(x_0, r)\times \sM$, where $x_0
 \in \R^d, r>0$, we have for any integer $K\in \N$, any $(x,i)\in B(x_0, r)\times \sM$,
 \begin{eqnarray}\label{iterationofu1}
u(x,i)&=&h(x,i)+\sum_{k=1}^{K}\sum_{\substack{l_1,...,l_k=1,\\l_1\neq i,l_2\neq l_1,...,l_k\neq l_{k-1}}}^m\wt{G}_{B(x_0, r)}^i(q_{il_1}(\wt{G}^{l_1}_{B(x_0, r)}q_{l_1l_2}(\dots (\wt{G}^{l_{k-1}}_{B(x_0, r)}q_{l_{k-1}l_k}h(\cdot,l_k)\dots)))(x)\nonumber\\
&&+\sum_{\substack{l_1,...,l_{K}=1,\\l_1\neq i,l_2\neq l_1,...,l_{K}\neq l_{K-1}}}^m \wt{G}_{B(x_0, r)}^i(q_{il_1}(\wt{G}^{l_1}_{B(x_0, r)} q_{l_1l_2}(\dots (\wt{G}^{l_{K-1}}_{B(x_0, r)} q_{l_{K-1}l_{K}}u(\cdot, l_{K}))\dots)))(x),\nonumber\\
\end{eqnarray}
where $h(x,i)$ is defined in \eqref{defofhi}.\\

When $\phi$ in \eqref{defintionofu} is bounded, for any $x_0\in \R^d, r\in (0, \delta_3]$, we have for any $(x,i)\in B(x_0, r)\times \sM$,
\begin{eqnarray}\label{iterationofu2}
u(x,i)&=&h(x,i)+\sum_{k=1}^{\infty}\sum_{\substack{l_1,...,l_k=1,\\l_1\neq i,l_2\neq l_1,...,l_k\neq l_{k-1}}}^m\wt{G}_{B(x_0, r)}^i(q_{il_1}(\wt{G}^{l_1}_{B(x_0, r)}q_{l_1l_2}(\dots (\wt{G}^{l_{k-1}}_{B(x_0, r)}q_{l_{k-1}l_k}h(\cdot,l_k))\dots)))(x),\nonumber\\
\end{eqnarray}
and we have 
  \begin{eqnarray}\label{representofuinstates}
u(x,i)&=&h(x,i)+\sum_{n=1}^{\infty}\sum_{l\in \sM}\sum_{(i,l_1,...,l_{n})\in \Psi(n;i,l)}\wt{G}_B^i(q_{il_1}(\wt{G}^{l_1}_{B} q_{l_1l_2}(\dots (\wt{G}^{l_{n-1}}_{B}q_{l_{n-1}l_{n}}h(\cdot,l_{n}))\dots)))(x),\nonumber\\
\end{eqnarray}
if we represent \eqref{iterationofu2} in terms of valid paths in each $\Psi(n;i,l)$. Here $\delta_3$ is the constant in Lemma \ref{norm0} and $\Psi(n;i,l)$ is defined in \eqref{path1}.
  \end{proposition} 
  \proof 
Define $\tau_2=\inf\{t> \tau_1:\Lambda_t\neq \Lambda_{\tau_1}\}$, where $\tau_1$ is defined in \eqref{defoftau1}.
  Then by Proposition \ref{representofu}, and the strong Markov property of $(X, \Lambda)$, we have for any $(x,i)\in B(x_0, r)\times\sM$,
  \begin{eqnarray*}\label{representationu2}
  u(x,i)&=&h(x,i)+\sum_{j=1,j\neq i}^m \wt{G}_{B(x_0, r)}^i(q_{ij}u(\cdot, j))(x)\nonumber\\
&=&h(x,i)+\sum_{j=1,j\neq i}^m\bE^{(x,i)}[\int_0^{\tau_{B(x_0 r)}\wedge \tau_1}q_{ij}u(X_s, j)ds]\nonumber\\
&=&h(x,i)+\sum_{j=1,j\neq i}^m \bE^{(x,i)}[\int_0^{\tau_{B(x_0 r)}\wedge \tau_1}q_{ij}(h_j(X_s)+\sum_{l=1,l\neq j}^m(\wt{G}_{B(x_0, r)}^jq_{jl}u(\cdot,l))(X_s))ds]\nonumber\\
&=&h(x,i)+\sum_{j=1,j\neq i}^m\wt{G}_{B(x_0, r)}^i(q_{ij}h(\cdot, j))(x)+\sum_{j,l=1, j\neq l,j\neq i}^m\wt{G}^i_{B(x_0, r)}(q_{ij}(\wt{G}^j_{B(x_0, r)}q_{jl}u(\cdot, l)))(x).\nonumber\\
   \end{eqnarray*}
 Iterating  \eqref{represenofu1} for $K-1$ times as above, we obtain \eqref{iterationofu1}.
 Then  by \eqref{GBQ}, if we represent \eqref{iterationofu1} in terms of vectors, we have for any $x\in B(x_0, r)$,
  \begin{eqnarray}\label{iterationofuvec}
u(x)=h(x)+\sum_{k=1}^{K}(\wt{G}_{B(x_0, r)}\wt{Q})^kh(x)+(\wt{G}_{B(x_0, r)}\wt{Q})^{K}u(x),
 \end{eqnarray} 
 where $u(x)=[u(x,1),...,u(x,m)]^T$, and $h(x)=[h(x,1),...,h(x,i)]^T.$
 
 If $\phi$ is bounded, then by \eqref{defintionofu}, $u$ is bounded. Then by Lemma \ref{norm0}, for any $x_0\in \R^d, r\in (0,\delta_3]$, any  bounded and nonnegative function $u$ on $B(x_0, r)$,
 $$\lim_{K\rightarrow\infty}(\wt{G}_{B(x_0, r)}\wt{Q})^{K}u(x)=\vec{0},$$ where $\vec{0}$ represents the zero vector of dimension $m$.
 Therefore,   
 \begin{equation*}\label{representofuinmatrix}
 u(x)=h(x)+\sum_{k=1}^{\infty}(\wt{G}_{B(x_0, r)}\wt{Q})^kh(x), \quad x\in B(x_0,r),
 \end{equation*}
i.e., for any $(x, i)\in B(x_0, r)\times \sM$,
\begin{eqnarray*}
u(x,i)=h(x,i)+\sum_{k=1}^\infty\sum_{\substack{l_1,...,l_n=1,\\l_1\neq i,l_2\neq l_1,...l_{k-1}\neq l_k}}^m \wt{G}_B^i(q_{il_1}(\wt{G}_B^{l_1}q_{l_1l_2}(...(\wt{G}_B^{l_{k-1}}q_{l_{k-1}l_n}h(\cdot, l_k))...)))(x).\nonumber\\
\end{eqnarray*}
 If we represent it in terms of valid paths, 
$$
u(x,i)=h(x,i)+\sum_{n=1}^{\infty}\sum_{l\in \sM}\sum_{(i,l_1,...,l_{n})\in \Psi(n;i,l)}^m \wt{G}_B^i(q_{il_1}(\wt{G}_B^{l_1}q_{l_1l_2}(...(\wt{G}_B^{l_{n-1}}q_{l_{n-1}l_n}h(\cdot, l_n))...)))(x).$$
\qed
\section {H\"older regularity} \label{S:3} 
\noindent{\it Proof of Theorem \ref{T:Holder}.} 

 For any $x_0\in \R^d, r\in (0, \wt{r}_0)$, denote $B=B(x_0, 2r/3)$, where $\wt{r}_0$ will be determined later. For any $i\in\sM$, let $\{X_s^{i},s\geq 0;\bP^x, x\in \R^d\}$ be the Hunt process corresponding to the operator $\sL_i$ starting from $x$. Define $\tau^{i}_B=\inf\{t>0:X_t^i\notin B \}.$
 
 First we claim that for any  $x\in B$, 
 \begin{equation*}
 h_i(x)=\wt{u}_i(x)+G^i_B(q_{ii}h_i)(x),
 \end{equation*}
 where $\wt{u}_i(x)=\bE^x[u(X_{\tau^i_B}^i,i)]$, $G_B^i$ is the Green operator for $\sL_i$ with zero Dirichlet boundary condition.
By the definition of harmonic functions with respect to $\sL_i+q_{ii}$,
\begin{eqnarray}\label{relationkilledandharmonic1}
 h_i(x)&=&\bE^{x}[\exp(\int_0^{\tau^{i}_B}q_{ii}(X^i_s)ds)u(X^i_{\tau^{i}_B},i)]\nonumber\\
 &=&\wt{u}_i(x)+\bE^{x}[(\exp(\int_0^{\tau^{i}_B}q_{ii}(X^i_s)ds)-1)u(X^i_{\tau^{i}_B},i)]\nonumber\\
 &=&\wt{u}_i(x)+\bE^{x}[\int_0^{\tau^{i}_B}(\exp(\int_s^{\tau^{i}_B}q_{ii}(X^i_t)dt))q_{ii}(X^i_s)u(X^i_{\tau^{i}_B},i)ds]\nonumber\\
  &=&\wt{u}_i(x)+\bE^{x}[\int_0^{\tau^{i}_B}q_{ii}(X^i_s)\bE^{X_s^i}[\exp(\int_0^{\tau^{i}_B}q_{ii}(X^i_t)dt)u(X^i_{\tau^{i}_B},i)]ds]\nonumber\\
 &=&\wt{u}_i(x)+\bE^{x}[\int_0^{\tau^{i}_B}q_{ii}(X^i_s)h_i(X^i_s)ds]\nonumber\\
 &=&\wt{u}_i(x)+\int_BG^i_B(x,y)q_{ii}(y)h_i(y)dy.\nonumber\\
  &=&\wt{u}_i(x)+G^i_Bq_{ii}h_i(x), 
 \end{eqnarray}
Notice that by Lemma \ref{supexittime}, for any $x\in B(x_0, r/2), r\in (0, r_0)$,
\begin{eqnarray}\label{e3.2}
|\int_BG^i_B(x,y)q_{ii}(y)h_i(y)dy|&\leq &4\wt{c}_1\Theta_5\|u\|_{\infty}r^2/9,
\end{eqnarray}
where $r_0$ is the constant in Lemma \ref{supexittime}. We will get the similar result when $x$ is replaced by $y$. 

Since by Theorem 1 in \cite{hinonlocal}, there exists a constant $C^{i}_3=C^{i}_3$, $r_4\in(0,1/4)$ and $\alpha\in (0,1)$ which both depend on$ (d,\Theta_1,\Theta_2,\Theta_3,\Theta_4,\beta)$ such that for any $x_0\in \R^d, r\in (0,r_4)$, any bounded function $\wt{u}_i$ defined in $\R^d$ that is harmonic with respect to $\sL_i$ in $B(x_0, r)$,
\begin{equation}\label{Holder1}
(\wt{u}_i(x)-\wt{u}_i(y))|\leq C^i_3\|u_i\|_{\infty}(\frac{\|x-y\|}{r})^{\alpha}, \hbox{for }x,y\in B(x_0, r/2).
\end{equation}
Let $\wt{r}_0=r_0\wedge r_4$.
 Therefore  by \eqref{relationkilledandharmonic1}, \eqref{e3.2}, \eqref{Holder1}, for any $r\in (0, \wt{r}_0)\subset (0,1/4), x,y\in B(x_0, r/2)$,
  \begin{eqnarray}
  |h_i(x)-h_i(y)|&\leq &|\wt{u}_i(x)-\wt{u}_i(y)|+|G^i_Bq_{ii}h_i(x)|+|G^i_Bq_{ii}h_i(y)|\nonumber\\
&\leq&C_3^i\|u\|_{\infty}(\frac{|x-y|}{r})^{\alpha}+8\wt{c}_1r^2\Theta_5\|u\|_{\infty}/9.
\end{eqnarray}
Set $r=|x-y|^{1/2}$, we have 
  \begin{eqnarray}\label{e3.5}
  |h_i(x)-h_i(y)|&\leq &|\wt{u}_i(x)-\wt{u}_i(y)|+|G^i_Bq_{ii}h_i(x)|+|G^i_Bq_{ii}h_i(y)|\nonumber\\
&\leq&C_3^i\|u\|_{\infty}|x-y|^{\alpha/2}+8\wt{c}_1|x-y|\Theta_5\|u\|_{\infty}/9\nonumber\\
&\leq &(C_3^i+8\Theta_5\wt{c}_1\wt{r}_0^{1-\alpha/2}/9)\|u\|_{\infty}|x-y|^{\alpha/2}.
\end{eqnarray}
For any $x_0\in \R^d, r\in (0,\wt{r}_0)$, denote $\wt{B}=B(x_0, 2r/3)$. Define $\tau_{\wt{B}}=\inf\{t>0:(X_t,\Lambda_t)\notin \wt{B}\times \sM\}.$ By Lemma \ref{supexittime}, for any $x \in 
B(x_0, r/2)$, we estimate each term in the second part of  \eqref{represenofu1} such that
\begin{eqnarray}\label{e3.6}
|\wt{G}^{i}_{\wt{B}}( q_{ik} u(\cdot,k))(x)|&\leq&4\wt{c}_1r^2\Theta_5\|u\|_{\infty}/9,
\end{eqnarray}
Therefore, by  \eqref{represenofu1}, \eqref{e3.5}, \eqref{e3.6}, for any $r\in (0, \wt{r}_0)$, $x,y \in B(x_0, r/2)$, we have
\begin{eqnarray}
\|u(x)-u(y)\|&=& \sup_{i\in \sM}|u(x,i)-u(y,i)|\nonumber\\
&\leq &\sup_{i\in \sM}(|h_i(x)-h_i(y)|+\sum_{k=1,k\neq i}^m(|G^{i}_{\wt{B}}( q_{ik} u(\cdot,k))(x)|+|G^{i}_{\wt{B}}( q_{ik} u(\cdot,k))(y)|)\nonumber\\
&\leq&(\max_{i\in\sM}C_3^i+8\Theta_5\wt{c}_1\wt{r}_0^{1-\alpha/2}/9)\|u\|_{\infty}|x-y|^{\alpha/2}+8(m-1)\wt{c}_1\Theta_5\|u\|_{\infty}r^2/9.\nonumber
\end{eqnarray}
Set $r=|x-y|^{1/2}$, $\alpha_1=\alpha/2.$ We have for any  $r\in (0,\wt{r}_0), x, y\in B(x_0, r/2)$,
\begin{eqnarray}
\|u(x)-u(y)\|&\leq& (\max_{i\in\sM}C_3^i+8\Theta_5\wt{c}_1\wt{r}_0^{1-\alpha/2}/9)\|u\|_{\infty}|x-y|^{\alpha/2}+8(m-1)\wt{c}_1\Theta_5\|u\|_{\infty}|x-y|/9\nonumber\\
&\leq&(\max_{i\in\sM}C_3^i+8m\wt{c}_1\wt{r}_0^{1-\alpha/2}/9)\Theta_5\|u\|_{\infty}|x-y|^{\alpha/2}:= C_1\|u\|_{\infty}|x-y|^{\alpha/2}.
\end{eqnarray}

\section {Scale-invariant Harnack inequality for each level}\label{S:4} 

 \begin{proposition}\label{killingEHI}
  Suppose that the operator $\mathcal{L} \in \sN(\Theta_1, \Theta_2, \Theta_3, \Theta_4, \gamma,\beta,c_1),$ and $\{\bar{X}_t,t\geq 0;\bP^x, x\in \R^d\}$ is the Hunt  process corresponding to the operator $\sL$ starting at $x$. Let $\wt{X}$ be the subprocess of $\bar{X}$  with killing rate $-q$ such that $\|q\|_{\infty}=\Theta_5$.
 There exists a constant $\wt{C}_2=\wt{C}_2(d, \beta, \gamma, \Theta_1, \vartheta)>0$ such that for any $x_0\in \R^d, r\in(0 , \delta_3]$, any nonnegative functions $h$ defined in $\R^d$ that is harmonic with respect to $\sL+q$ in $B(x_0, r)$, 
\begin{equation}\label{killedEHI}
h(y_1)\leq C_2h(y_2), \quad \hbox{for any }y_1, y_2\in B(x_0, r/4).
\end{equation}
Here $\delta_3$ is the constant in Lemma \ref{norm0}.
 \end{proposition}
\proof For simplicity, we just give a proof for the $d\geq 3$, the case for $d=2$ is similar by using Lemma \ref{2.9} instead of Lemma \ref{2.19}. This follows the same idea of \cite[Theorem 1.7]{BHP}. We spell out the details here for the reader's convenience.

For any $x_0\in \R^d, r\in (0, \delta_3]$,where $\delta_3$ is the constant in Lemma \ref{norm0},
let $\tau_{B(x_0, r)}=\inf\{t> 0: \bar{X}_t\notin B(x_0, r)\}$.
 By shrinking the value of $r$ a little bit, we may assume without loss of generality that 
$h(x) = \bE^x  \big[ h(\wt{X}_{\tau_{B(x_0,r)}}) \big], x\in B(x_0, r) $.
As $h$ is the increasing limit of bounded harmonic functions $h_n (x):= \bE^x  \big[ (h\wedge n)(\wt{X}_{\tau_{B(x_0, r)}}) \big]$,
it suffices to establish the Harnack inequality for bounded nonnegative harmonic functions $h$ in $B(x_0, r)$. 
 
Denote $B_r=B(x_0, r)$. Then by Theorem \ref{kiillingharmonicrep},  there exists a unique measure $\mu_h$ on $\partial B_r$, such that for any for any $x\in B(x_0, r)$, 
 \begin{eqnarray}\label{rightharmonic}
  h(x)&=&\int_{\partial B_r}\wt{M}_{B_r}(x,z)h(z)\mu_{h}(dz)+\int_{\overline{B_r}^c}\int_{B_r}\wt{G}_{ B_r}(x,y)J^{b_2}(y,z)dyh(z)dz\nonumber\\
  &=&\int_{\partial B_r}\wt{M}_{B_r}(x,z)h(z)\mu_{h}(dz)+\int_{\overline{B_r}^c}\int_{B_r\setminus B_{r/2}}\wt{G}_{ B_r}(x,y)J^{b_2}(y,z)dyh(z)dz\nonumber\\
  &&+\int_{\overline{B_r}^c}\int_{ B_{r/2}}\wt{G}_{B_r}(x,y)J^{b_2}(y,z)dyh(z)dz.\nonumber\\
  &=&h_1(x)+h_2(x)+h_3(x).
  \end{eqnarray}
   
  Since for any  $z\in \partial B_r$, $M^{\Delta}_{B_r}(\cdot, z)$ is harmonic with respect to $\Delta$ in $B(x_0, r)$, then there exists a constant $C_2=C_2(d)>0$ such that for any $x,x'\in B(x_0,r/2)$,
  \begin{eqnarray}\label{h1}
   M^{\Delta}_{B_r}(x,z)\leq C_2M^{\Delta}_{B_r}(x',z).
   \end{eqnarray}
Therefore, by Lemma \ref{killmartinkernel}, 

\begin{equation}\label{h11}
h_1(x)\leq C_0^4C_2h_1(x')\hbox{ for any }x,x'\in B(x_0,r/2).
\end{equation}

 Also, by Lemma \ref{2.19}, there exists a constant $C_3=C_3(d)>0$ such that 
  \begin{eqnarray}\label{h2}
   G^{\Delta}_{B_r}(x,y)\leq C_3G^{\Delta}_{B_r}(x', y) \hbox{ for any }x,x'\in B(x_0,r/4), y\in B_r\setminus B_{r/2}.
     \end{eqnarray}
   Therefore, by Lemma \ref{killgreenfunction}, \eqref{h2}, $h_2(x)\leq C_0^2C_3h_2(x')$ for any $x,x'\in B(x_0,r/4)$.\\
Next, on one hand, by the UJS condition \eqref{UJScondition} for $J^{b_2}$, \eqref{e:3.45}, Lemma \ref{killgreenfunction}, for any $x\in B(x_0, r/2)$,
\begin{eqnarray}\label{h31}
h_3(x)&=&\int_{\overline{B_r}^c}\int_{ B_{r/2}}\wt{G}_{B_r}(x,y)J^{b_2}(y,z)dyh(z)dz\nonumber\\
&\leq &C_0\vartheta\int_{\overline{B_r}^c}\int_{ B_{r/2}}G^{\Delta}_{B_r}(x,y)\frac{1}{|B(y,r/4)|}\int_{B(y,r/4)}J^{b_2}(u,z)dudyh(z)dz \nonumber\\
&\leq&4^dC_0\vartheta/(\omega_d r^d)\int_{\overline{B_r}^c}(\int_{ B_{r/2}}G^{\Delta}_{B_r}(x,y)dy)\int_{B_{3r/4}}J^{b_2}(u,z)duh(z)dz\nonumber\\
&\leq&c_3C_0K_1\vartheta r^{2-d}\int_{\overline{B_r}^c}\int_{B_{3r/4}}J^{b_2}(u,z)duh(z)dz,
\end{eqnarray}
where $\omega_d$ and $c_3$ are constants depending only on $d$, and therefore $c_4C_0K_1\vartheta /\omega_d$ depends on $(d,\beta, \gamma, \Theta_1,\vartheta)$, and  independent of $r$.\\
On the other hand, by Lemma \ref{2.19}, Lemma \ref{killgreenfunction}, for any $x'\in B(x_0, r/2)$,
\begin{eqnarray}\label{h32}
h(x')&\geq &\bE^{x'}[h(\wt{X}_{\tau_{B_{7r/8}}}); \wt{X}_{\tau_{B_{7r/8}}}\in \overline{B_r}^c]\nonumber\\
&\geq&C_0^{-1}\int_{\overline{B_r}^c}\Big(\int_{B_{3r/4}}G_{B_{7r/8}}^{\Delta}(x^{'},y)J^{b_2}(y,z)dy\Big)h(z)dz\nonumber\\
&\geq &c_4(C_0K_1)^{-1}r^{2-d}\int_{\overline{B_r}^c}\int_{B_{3r/4}}J^{b_2}(y,z)dyh(z)dz,
\end{eqnarray}
where $c_4(C_0K_1)^{-1}$  depends only on $(d, \gamma, \Theta_1)$ and  independent of $r$.\\
Therefore by \eqref{rightharmonic}, \eqref{h1},\eqref{h2},\eqref{h31},\eqref{h32}, for any $x,x' \in B(x_0, r/4)$,
$$h_1(x)+h_2(x)+h_3(x)\leq (C_2+C_3+c_3C_0^2K_1^2\vartheta/c_4)h(x'):=\wt{C}_2h(x'),$$
where $\wt{C}_2$ depends only on $(d,\beta,\gamma,\Theta_1,\vartheta).$
 
  \qed
 
  \medskip
\noindent{\it Proof of Theorem \ref{T:1}}:\\ 
Fix $x_0\in \R^d, r\in (0, \wt{r}_1]$, with $\wt{r}_1=r_0\wedge \delta_3$, where $r_0$ and $\delta_3$ are the constants in Lemma \ref{supexittime} and Lemma \ref{norm0} respectively. Notice that $\wt{r}_1< r_5 <1$ by Lemma \ref{estimationforexp} and \eqref{choiceofr}.

Similar to the beginning of the argument in Proposition \ref{killingEHI}, it suffices to establish the Harnack inequality for bounded nonnegative functions $u$ defined in $\R^d$ that is $\mathcal{G}$-harmonic  in $B(x_0, r)$. For any $k\in \sM$, let $\wt{X}^k$ be the Hunt process corresponding to the operator $\sL_k+q_{kk}$ for $\sL_k \in \sN(\Theta_1, \Theta_2, \Theta_3, \Theta_4, \gamma,\beta,c_1)$.

For any $x_0\in \R^d, r\in (0,  \wt{r}_1]$, denoting $B=B(x_0,r/2)$,  by the definition of $h(x,k)$ and $\tau_1$  in \eqref{defofh} and \eqref{defoftau1},  \eqref{MP1}, $h(x,k)$ is harmonic with respect to $\sL_k+q_{kk}$ in $B$, i.e, 
$$h(x,k)=\bE^{(x,k)}[u(\wt{X}^k_{\tau_{B}},k)] \hbox{ for }x\in B.$$
Then by Theorem \ref{kiillingharmonicrep}, there exists a unique measure $\mu_{u}^{k}$ on $\partial B$, such that for any $x\in B$, 
$$h(x,k)=\int_{\overline{B}^c}\int_{B}
\wt{G}^k_{B}(x,y)J^{b_2,k}(y, z)dyu(z,k)dz+\int_{\partial B}\wt{M}^k_{B}(x,z)u(z, k)\mu^k_u(dz).$$
Hence by Proposition \ref{killingEHI}, there exists a constant $\wt{C}_2=\wt{C}_2(d, \beta, \gamma, \Theta_1, \vartheta)>0$
 such that for any $x,x'\in B(x_0, r/8)$, 
  \begin{eqnarray}\label{killingEHI1}
 h(x,k)\leq \wt{C}_2h(x',k).
  \end{eqnarray}
For any $C>0, x\in B,k,l\in \sM, r\in (0, \wt{r}_1]$, define
 \begin{eqnarray}\label{Fform1}
 F_{kl}^{I}(C;x,z,r):&=&\sum_{n=1}^{\infty}\sum_{\Psi^0(n;k,l)}C^nq^0_{kl_1}q^0_{l_1l_2}...q^0_{l_{n-1}l}\int_{B^{n}}G^{\Delta}_B(x,y_1)G^{\Delta}_B(y_1,y_2)...G^{\Delta}_B(y_{n-1},y_n) (\int_{B}\nonumber\\
&&CG^{\Delta}_{B}(y_n,y)J^{b_2,l}(y, z)\mathds{1}_{\overline{B}^c}(z)dy)dy_1dy_2...dy_n, \hbox{ for }z\in \overline{B}^c,
\end{eqnarray}
and 
 \begin{eqnarray}\label{Fform2}
F_{kl}^{II}(C;x,z,r):&=&\sum_{n=1}^{\infty}\sum_{\Psi^0(n;k,l)}C^nq^0_{kl_1}q^0_{l_1l_2}...q^0_{l_{n-1}l}\int_{B^n}G^{\Delta}_B(x,y_1)G^{\Delta}_B(y_1,y_2)...G^{\Delta}_B(y_{n-1},y_n)\nonumber\\
&&(C^2M^{\Delta}_{B}(y_n,z)\mathds{1}_{\partial B}(z))dy_1dy_2...dy_n, \hbox{ for }z\in \partial B.
\end{eqnarray}
For each $k\in \sM$, recall that
\begin{equation}\label{ek}
E(k)=\{l\in \sM\setminus k:\Psi^0(n,k,l)\neq \emptyset,n\in \N\}.
\end {equation}
Then by \eqref{condQ}, Lemma \ref{killgreenfunction}, Lemma \ref{killmartinkernel},  Theorem \ref{kiillingharmonicrep}, \eqref{representofuinstates}, \eqref{Fform1},\eqref{Fform2},\eqref{ek}, for any $(x,k)\in B\times \sM$, it yields that
\begin{eqnarray}\label{ucomparision}
&&h(x,k)+ \sum_{l\in E(k)\cup \{k\}}(\int_{\overline{B}^c}F^{I}_{kl}(c_0C_0^{-1};x,z,r)u(z,l)dz+\int_{\partial B}F^{II}_{kl}(c_0C_0^{-1};x,z,r)u(z,l)\mu^l_u(dz))\nonumber\\
&\leq&u(x,k)\leq h(x,k)+ \sum_{l\in E(k)\cup\{k\}}(\int_{\overline{B}^c}F^I_{kl}(C_0;x,z,r)u(z,l)dz+\int_{\partial B}F^{II}_{kl}(C_0;x,z,r)u(z,l)\mu^l_u(dz)).\nonumber\\
\end{eqnarray}
By the 3G-Lemma (Lemma \ref{lemma3.14}), for any  $x,y\in B, $
 \begin{eqnarray}\label{e:4.13}
 \int_B G_B^{\Delta}(x,y_1)G_B^{\Delta}(y_1,y)dy_1&\leq& C_GG_B^{\Delta}(x,y)(\sup_{x\in B}\int_{B}|x-y_1|^{2-d}dy_1+\sup_{y\in B}\int_{B}|y-y_1|^{2-d}dy_1)\nonumber\\
 &\leq &C_G\omega_dr^2 G_B^{\Delta}(x,y)=:C_gr^2G_B^{\Delta}(x,y),
 \end{eqnarray} 
 Then by Lemma \ref{lemma3.14},
  \eqref{e:4.13}, for any $r\in (0, \wt{r}_1], x\in B, z\in \overline{B}^c$, we have
 \begin{eqnarray}\label{Fformup}
  F_{kl}^I(C_0;x,z,r)&\leq &\sum_{n=1}^{\infty}\sum_{\Psi^0(n;k,l)}(C_gC_0)^{n}q^0_{kl_1}q^0_{l_1l_2}...q^0_{l_{n-1}l}r^{2n}(\int_{B}C_0G^{\Delta}_B(x,y)J^{b_2,l}(y,z)\mathds{1}_{\overline{B}^c}(z)dy) \nonumber\\
  &\leq &\sum_{n=1}^{\infty}(C_gC_0r^2)^{n}\sum_{\Psi^0(n;k,l)}q^0_{kl_1}q^0_{l_1l_2}...q^0_{l_{n-1}l}(\int_{B}C_0G^{\Delta}_B(x,y)J^{b_2,l}(y,z)\mathds{1}_{\overline{B}^c}(z)dy).\nonumber\\
    \end{eqnarray}
 
 By Lemma \ref{2.19},
for any $x,y\in B(x_0, 3r/8), B=B(x_0, r/2),$
\begin{eqnarray}\label{greenfuncestm1}
G^{\Delta}_B(x,y)\geq K_1^{-1}(3r/4)^{2-d}/36,
\end{eqnarray}
and 
for any  $x\in B(x_0, r/8), y\in B(x_0, 3r/8),$
\begin{eqnarray}\label{greenfuncestm2}
G^{\Delta}_B(x,y)\geq 3K_1^{-1}(r/2)^{2-d}/16\geq K_1^{-1} (3r/4)^{2-d}/36.
\end{eqnarray}
Then by \eqref{greenfuncestm1}, \eqref{greenfuncestm2}, for any $k,l\in \sM, r\in (0, \wt{r}_1], x\in B, z\in \overline{B}^c$, we have
      \begin{eqnarray}\label{Fkl1lowbd}
  &&F_{kl}^{I}(c_0C_0^{-1};x,z,r)\nonumber\\
  &\geq& \sum_{n=1}^{\infty}\sum_{\Psi^0(n;k,l)}c_0^nC_0^{-n}q^0_{kl_1}q^0_{l_1l_2}...q^0_{l_{n-1}l}\int_{B(x_0,3r/8)^{n-1}}G_B^{\Delta}(x,y_1)...\int_{B(x_0,r/8)}G_B^{\Delta}(y_{n-1},y_n)(\int_{B}C_0^{-1}G_B^{\Delta}(y_n,y)\nonumber\\
  &&J^{b_2,l}(y,z)\mathds{1}_{\overline{B}^c}(z)dy)dy_ndy_{n-1}...dy_1\nonumber\\
&\geq &\sum_{n=1}^{\infty}\sum_{\Psi^0(n;k,l)}c_0^nC_0^{-n}q^0_{kl_1}q^0_{l_1l_2}...q^0_{l_{n-1}l}(K_1^{-1}\int_{B(x_0,3r/8)}((3r/4)^{2-d}/36) dy_1)^{n-1}\int_{B(x_0, r/8)}(K_1^{-1}(3r/4)^{2-d}/36)\nonumber\\
 &&(\int_{B}C_0^{-1}G_B^{\Delta}(y_n,y)J^{b_2,l}(y,z)\mathds{1}_{\overline{B}^c}(z)dy)dy_n\nonumber\\
   &\geq &\sum_{n=1}^{\infty}\sum_{\Psi^0(n;k,l)}c_0^nC_0^{-n}q^0_{kl_1}q^0_{l_1l_2}...q^0_{l_{n-1}l}((36K_1)^{-1}(3r/4)^{2-d}\cdot \omega_d(r/8)^d)^{n}(\inf_{y_n\in B(x_0,r/8)}\int_{B}C_0^{-1}G_B^{\Delta}(y_n,y)\nonumber\\
  && J^{b_2,l}(y,z)\mathds{1}_{\overline{B}^c}(z)dy)\nonumber\\
&\geq&\sum_{n=1}^{\infty}(c_0K_dr^2/(K_1C_0))^{n}\sum_{\Psi^0(n;k,l)} q^0_{kl_1}q^0_{l_1l_2}...q^0_{l_{n-1}l}(\inf_{y_n\in B(x_0,r/8)}\int_{B}C_0^{-1}G_B^{\Delta}(y_n,y)J^{b_2,l}(y,z)\mathds{1}_{\overline{B}^c}(z)dy),\nonumber\\
   \end{eqnarray}
   where $K_d=K_d(d)$ is some constant  such that $K_d\leq C_g$ in \eqref{e:4.13}.
   
      For any $s\geq 0, k,l\in\sM$, define 
      \begin{eqnarray}\label{Hform}
    H_{kl}(s):= \sum_{n=1}^{\infty}a_n(kl)s^{n},
        \end{eqnarray}
        where $ a_n(kl)=\sum_{\Psi^0(n;k,l)}q^0_{kl_1}q^0_{l_1l_2}...q^0_{l_{n-1}l}.$
        Since each $a_n\leq \Theta_5^n(m-1)^n$,  then by the choice of $\delta_3$ in \eqref{choiceofr},   for any $r\in (0,\wt{r}_1]\subset (0, \delta_3)$,
\begin{equation}\label{Hklbound}
H_{kl}(s)\leq \sum_{n=1}^{\infty}((m-1)\Theta_5s)^n< 2 \hbox { for any }s \in (0, C_gC_0\wt{r}_1^2].
\end{equation}

              Then by \eqref{Fformup},  \eqref{Hform},  for any $ r\in (0, \wt{r}_1]$, $x\in B,z\in \overline{B}^c,k,l\in \sM$,
        \begin{eqnarray}\label{Fkl1up}
         F_{kl}^I(C_0;x,z,r)&\leq& H_{kl}(C_gC_0r^2)(C_0\int_{B}G^{\Delta}_B(x,y)J^{b_2,l}(y,z)\mathds{1}_{\overline{B}^c}(z)dy).
             \end{eqnarray}
        Similarly, by \eqref{martinkernel}, \eqref{Fform2},\eqref{Fformup}, \eqref{Fkl1up},
         for any $ r\in (0, \wt{r}_1], x\in B, z\in \partial B, k,l \in \sM,$ 
           \begin{eqnarray}\label{Fkl2up}
         F_{kl}^{II}(C_0;x,z,r) &\leq &H_{kl}(C_gC_0r^2)(C_0^2M^{\Delta}_{B}(x,z)\mathds{1}_{\partial B}(z)).
         \end{eqnarray}       
   Also, by \eqref{Fkl1lowbd},  \eqref{Hform}, for any $ r\in (0, \wt{r}_1]$, $x\in B,z\in \overline{B}^c,k,l\in \sM$,
   \begin{eqnarray}\label{Fkl1lowbd2}
   F_{kl}^{I}(c_0C_0^{-1};x,z,r)&\geq &H_{kl}(c_0K_dr^2/(K_1C_0))(\inf_{\wt{y}\in B(x_0,r/8)}\int_{B}C_0^{-1}G_B^{\Delta}(\wt{y},y)J^{b_2,l}(y,z)\mathds{1}_{\overline{B}^c}(z)dy)\nonumber\\
   \end{eqnarray}
Similarly, by \eqref{Fform2}, \eqref{Fkl1lowbd} \eqref{Fkl1lowbd2}, for any $ r\in (0, \wt{r}_1], x\in B, z\in \partial B,k\in \sM,l\in E(k)$,
\begin{equation}\label{Fk2lowbd}
F_{kl}^{II}(c_0C_0^{-1};x,z,r)\geq H_{kl}(c_0K_dr^2/(K_1C_0))(C_0^{-2}\inf_{\wt{y}\in B(x_0,r/4)}M_B^{\Delta}(\wt{y},z)\mathds{1}_{\partial{B}}(z)).
\end{equation}
Recall that  $m_{k l}:=\inf\{n\geq 0: \Psi^0(n,k,l)\neq \emptyset\}$, the smallest integer $n$ such that one can go from level $k$ to level $l$ in $n$ steps.
        By \eqref{Hklbound} , for any $s \in (0, C_gC_0\wt{r}_1^2]$,    
       \begin{equation}\label{estofHkl}
       a_{m_{kl}}s^{m_{kl}} \leq  H_{kl}(s)=a_{m_{kl}}s^{m_{kl}}(1+o(s))\leq 3a_{m_{kl}}s^{m_{kl}},
       \end{equation}
       where $o(s)$ stands for some function of $s$ converging to $0$  as $s$ approaches $0$.

Hence by \eqref{killingEHI1},\eqref{ucomparision},\eqref{Hklbound}, \eqref{Fkl1up},\eqref{Fkl2up}, \eqref{estofHkl}, Lemma \ref{killgreenfunction}, Lemma \ref{killmartinkernel}, Theorem \ref{kiillingharmonicrep}, for any $x,x'\in B(x_0, r/8), k\in \sM$,
\begin{eqnarray}\label{upbd}
u(x,k)&\leq&h(x,k)+ \sum_{l\in E(k)\cup\{k\}}H_{kl}(C_gC_0r^2)(C_0\int_{\overline{B}^c}\int_{B}G^{\Delta}_B(x,y)J^{b_2,l}(y,z)u(z,l)dydz\nonumber\\
&&+\int_{\partial B} C_0^2M^{\Delta}_{B}(x,z)u(z,l)\mu^l_u(dz))\nonumber\\
&\leq&(1+ 3a_{m_{kl}}(C_gC_0r^2)^{m_{kl}}C_0^4)h(x,k)+ \sum_{l\in E(k)} 3a_{m_{kl}}(C_gC_0r^2)^{m_{kl}}C_0^4h(x,l)\nonumber\\
&\leq&\wt{C}_2 (1+3a_{m_{kl}}(C_gC_0r^2)^{m_{kl}}C_0^4)h(x',k)+ \sum_{l\in E(k)} 3\wt{C}_2C_0^4a_{m_{kl}}(C_gC_0)^{m_{kl}}r^{2m_{kl}}h(x',l).\nonumber\\
\end{eqnarray}
Let $C_3=\max_{k\in \sM}(\sum_{l\in E(k)\cup \{k\}} 3\wt{C}_2C_0^4a_{m_{kl}}(C_gC_0)^{m_{kl}}+\wt{C}_2).$
Since each $m_{kl}$ is finite, then $C_3$ is finite. By \eqref{upbd}, for any $x,x'\in B(x_0, r/8), k\in \sM$,
\begin{eqnarray}\label{upbd1}
u(x,k)\leq C_3(h(x',k)+\sum_{l\in E(k)}r^{2m_{kl}}h(x',l)).
\end{eqnarray}
 Since by \eqref{killingEHI1} and Theorem \ref{kiillingharmonicrep}, for any $k\in\sM$, $x' \in B(x_0, r/8)$, 
\begin{eqnarray}\label{4.27}
&& \wt{C}_2 \inf_{\wt{y}\in B(x_0, r/8)}\Big(\int_{\overline{B}^c}\int_{B} \wt{G}^k_{B}(\wt{y},y)J^{b_2,k}(y, z)dyu(z,k)dz+\int_{\partial B}\wt{M}^k_{B}(\tilde{y},z)u(z, k)\mu^k_u(dz)\Big),\nonumber\\
&\geq &\int_{\overline{B}^c}\int_{B}\wt{G}^k_{B}(x',y)J^{b_2,k}(y, z)dyu(z,k)dz+\int_{\partial B}\wt{M}^k_{B}(x',z)u(z, k)\mu^k_u(dz),\nonumber\\
\end{eqnarray}
then by \eqref{killingEHI1}, \eqref{Fkl1lowbd2}, \eqref{Fk2lowbd},\eqref{estofHkl},\eqref{4.27}, Lemma \ref{killgreenfunction}, Lemma \ref{killmartinkernel}, Theorem \ref{kiillingharmonicrep}, for any $x,x'\in B(x_0, r/8),k\in \sM$,
\begin{eqnarray}\label{lowbd}
u(x,k)&\geq&h(x,k)+ \sum_{l\in E(k)\cup\{k\}}H_{kl}(c_0K_dr^2/(K_1C_0))(C_0^{-1}\int_{\overline{B}^c}\int_{B}G^{\Delta}_B(x,y)J^{b_2,l}(y,z)u(z,l)dydz\nonumber\\
&&+\int_{\partial B} C_0^{-2}M^{\Delta}_{B}(x,z)u(z,l)\mu^l_u(dz))\nonumber\\
&\geq& h(x,k)+\sum_{l\in E(k)} a_{m_{kl}}(c_0K_dr^2/(K_1C_0))^{m_{kl}}C_0^{-4}h(x,l)\nonumber\\
&\geq& h(x,k)+\sum_{l\in E(k)} a_{m_{kl}}(c_0K_d/(K_1C_0))^{m_{kl}}r^{2m_{kl}}C_0^{-4} h(x,l)\nonumber\\
&\geq&\wt{C}_2h(x',k)+ \sum_{l\in E(k)}(\wt{C}_2C_0^4)^{-1}a_{m_{kl}}(c_0K_d/(K_1C_0))^{m_{kl}}r^{2m_{kl}}h(x',l).
\end{eqnarray}
Let  $C_2=\min\{(\wt{C}_2C_0^4)^{-1}a_{m_{kl}}(c_0K_d/(K_1C_0))^{m_{kl}},l\in E(k),k\in \sM\} \wedge (1/\wt{C}_2),$ which is also finite.
Then by \eqref{lowbd}, for any $x,x'\in B(x_0, r/8),k\in \sM$,
\begin{eqnarray}\label{lowbd1}
u(x,k)\geq C_2(h(x',k)+\sum_{l\in E(k)}r^{2m_{kl}}h(x',l)).
\end{eqnarray}
Here both $C_2$ and $C_3$ depend on $(d,m, c_0,\beta, \gamma, \Theta_1,\Theta_2,\Theta_3,\Theta_4,\vartheta ,Q^0)$, and independent of $r$.

In particular, by \eqref{upbd1},\eqref{lowbd1}, for any $x_0\in \R^d, r\in (0,\wt{r}_1], x,x'\in B(x_0, r/8),k\in \sM$,
\begin{eqnarray}\label{onelevel}
u(x,k)\leq (C_3/C_2)C_2(h(x',k)+\sum_{l\in E(k)}r^{2m_{kl}}h(x',l))\leq (C_3/C_2)u(x',k)=:Cu(x',k).
\end{eqnarray}
\qed
 \section {Scale-invariant Full Rank Harnack inequality}\label{S:4} 
\noindent{\it Proof of Theorem \ref{T:2}.}

Fix $x_0\in \R^d$ and $r\in (0, \wt{r}_1)$,where $\wt{r}_1$ is the constant in Theorem \ref{T:1}. Similar to the beginning of the argument in Proposition \ref{killingEHI}, it suffices to establish the Harnack inequality for bounded nonnegative functions $u$ defined in $\R^d$ that is $\mathcal{G}$-harmonic  in $B(x_0, r)$.
Since $Q$ is strictly irreducible, then by \eqref{condQ}, \eqref{lowerbdofQ0},  we know that 
\begin{eqnarray}\label{qbound}
0<q_0\leq q_{kl}^0 , \hbox { for any } k,l\in \sM, k\neq l.
\end{eqnarray}

Also,  in \eqref{upbd1} \eqref{lowbd}, $E(k)=\sM$ for each $k\in \sM$, such that
\begin{equation}\label{eqmkl1}
m_{kl}=1, a_{m_{kl}}=q^0_{kl}\hbox{ for each }l\in\sM\setminus \{k\},\hbox{ and }m_{kk}=2.
\end{equation}
Then by 
\eqref{qbound},\eqref{eqmkl1}, for any $x,x'\in B(x_0, r/8) , k\in \sM$, we have 
\begin{eqnarray}\label{lowerbdu}
u(x,k)\geq (\wt{C}_2C_0^4)^{-1}c_0q_0K_dr^2/(K_1C_0)\sum_{l\in\sM\setminus \{k\}}h(x',l):=K_3r^2\sum_{l\in\sM\setminus \{k\}}h(x',l),
\end{eqnarray}
where $K_3=(\wt{C}_2C_0^4)^{-1}c_0q_0K_dr^2/(K_1C_0)$.
Since $\wt{r}_1<1$ by Theorem \ref{T:1}, then by \eqref{upbd1}, \eqref{qbound},\eqref{eqmkl1}, for any $r\in (0, \wt{r}_1], x,x'\in B(x_0, r/8), k\in \sM$, we have 
\begin{eqnarray}\label{upperbdu}
u(x,k)&\leq& C_3h(x',k)+ C_3r^2\sum_{l\in\sM\setminus \{k\}}  h(x',l).
\end{eqnarray}
Hence by Theorem \ref{T:1}, \eqref{lowerbdu}, \eqref{upperbdu}, for any $r\in (0, \wt{r}_1], x,x'\in B(x_0, r/8)$, any $k,i\in \sM, k\neq i$,
\begin{eqnarray}\label{e:5.6}
u(x,k)&\asymp& u(x',k)\leq C_3r^2 \sum_{l\neq k,i} h(x',l)+C_3r^2h(x',i)+C_3h(x',k)\nonumber\\
&\leq&(C_3/K_3)K_3r^2\sum_{l\neq k,i} h(x',l)+C_3u(x',i)+C_3h(x',k)\nonumber\\
&\leq &(C_3/K_3+2C_3)(u(x',i)+h(x',k))=\wt{C}_1(u(x',i)+h(x',k)),
\end{eqnarray}
where $\wt{C}_1:=C_3/K_3+2C_3+C$ depending only on $(d,\beta,\gamma, m, c_0,q_0,\Theta_1,\Theta_2,\Theta_3,\Theta_4,\vartheta,Q^0)$ and independent of $r$.

Furthermore, by  \eqref{lowerbdu},\eqref{upperbdu}, for any $ r\in (0,\wt{r}_1], x,x'\in B(x_0, r/8)$, any $k,i\in \sM,k\neq i$,
\begin{eqnarray*}\label{e:5.8}
u(x,k)&\leq &C_3(1+\wt{r}_1^2)\sum_{l\in \sM}h(x',l) \leq C_3(1+\wt{r}_1^2)(r^{-2}h(x',i)+\sum_{l\in\sM\setminus\{i\}}h(x', l))\nonumber\\
&\leq& C_3(1+\wt{r}_1^2)(r^{-2}u(x',i)+(1/K_3)r^{-2}u(x',i))\leq C_3(1+\wt{r}_1^2)(1+1/K_3)r^{-2}u(x',i)\nonumber\\
&\leq &\wt{C}_3(1+\wt{r}_1^2)r^{-2}u(x',i).
\end{eqnarray*}
\qed

\end{document}